\documentclass[12pt,a4paper]{article}

\usepackage[left=1cm,right=1cm,top=2cm,bottom=2cm]{geometry}
\usepackage[dvips]{graphics}
\usepackage[dvips]{epsfig}
\usepackage{color}
\usepackage{hyperref}
\usepackage{tikz,authblk}

\newtheorem{lemma}{Lemma}
\newtheorem{thm}[lemma]{Theorem}

\newtheorem{prop}[lemma]{Proposition}
\newtheorem{defn}{Definition}
\newtheorem{rem}{Remark}

\newtheorem{example}{Example}

\newcommand{\dimo}{\noindent \emph{Proof. }}
\newcommand{\qed}{\\ \rightline{$\Box$ \ \ \ \ \ \ \ \ \ \ \ \ \ \ \ }\\}

\newcommand{\G}{\Gamma}
\newcommand{\g}{\gamma}

\usepackage{latexsym}
\usepackage{amsfonts}
\usepackage{amssymb}
\usepackage{amsmath}
\usepackage[english]{babel}

\begin{document}

\title{Kirby diagrams and 5-colored graphs representing compact 4-manifolds}

 \renewcommand{\Authfont}{\scshape\small}
 \renewcommand{\Affilfont}{\itshape\small}
 \renewcommand{\Authand}{ and }

\author[1] {Maria Rita Casali}
\author[2] {Paola Cristofori}

\affil[1] {Department of Physics, Mathematics and Computer Science, University of Modena and Reggio Emilia \ \ \ \ \ \ \ \ \ \ \   Via Campi 213 B, I-41125 Modena (Italy), casali@unimore.it}

\affil[2] {Department of Physics, Mathematics and Computer Science, University of Modena and Reggio Emilia  \ \ \ \ \ \ \ \ \ \ \   Via Campi 213 B, I-41125 Modena (Italy), paola.cristofori@unimore.it}

\maketitle

\abstract{It is well-known that in dimension 4 any framed link $(L,c)$ uniquely represents the PL 4-manifold $M^4(L,c)$ obtained  from $\mathbb D^4$ by adding 2-handles along $(L,c)$. 
Moreover, if trivial dotted components are also allowed (i.e. in case of a {\it Kirby diagram} $(L^{(*)},d)$), the associated PL 4-manifold $M^4(L^{(*)},d)$ is obtained from $\mathbb D^4$ by adding 1-handles along the dotted components and 2-handles along the framed components. \\
In this paper  we study the relationships between framed links and/or Kirby diagrams and the representation theory of compact PL manifolds by edge-colored graphs: in particular, we describe how to construct algorithmically a (regular) 5-colored graph representing $M^4(L^{(*)},d)$, directly ``drawn over" a planar diagram of $(L^{(*)},d)$,  or equivalently how to algorithmically obtain a triangulation of $M^4(L^{(*)},d)$.  
As a consequence,  the procedure yields triangulations for any closed (simply-connected) PL 4-manifold admitting handle decompositions without 3-handles.   
\\  
Furthermore, upper bounds for both the invariants gem-complexity and regular genus of  $M^4(L^{(*)},d)$ are obtained, in terms of the combinatorial properties of the Kirby diagram.}
\endabstract
 


\maketitle

\section{Introduction}

Among combinatorial tools representing PL manifolds,  {\it framed links} (and/or  {\it Kirby diagrams}) turn out to be a very synthetic one, both in the 3-dimensional setting 
and in the 4-dimensional one, while {\it edge-colored graphs} have the advantage to represent all compact PL manifolds and to allow the definition and computation of 
interesting PL invariants in arbitrary dimension (such as the {\it regular genus}, which extends the Heegard genus, and the {\it gem-complexity}, similar to Matveev's complexity of a 3-manifold). 

\smallskip

Previous works exist establishing a connection between the two theories, both in the 3-dimensional and 4-dimensional setting (\cite{Lins-book}, \cite{Casali JKTR2000}, \cite{Casali Compl2004}): they make use of the so called {\it edge-colored graphs with boundary}, which are dual to colored triangulations of PL manifolds with non-empty boundary, and fail to be regular. More recently, a unifying method has been introduced and studied, so to represent by means of regular colored graphs all compact PL manifolds, via the notion of {\it singular manifold} associated to a PL manifold with non-empty boundary. 
  
Purpose of the present work is to update the relationship between framed links/Kirby diagrams and colored graphs (or, equivalently, colored triangulations) in dimension 4, by making use of regular 5-colored graphs representing compact PL 4-manifolds. The new tool turns out to be significantly more efficient than the classic one, both with regard to the simplicity and algorithmicity of the procedure and with regard to the possibility of estimating graph-defined PL invariants directly from the Kirby diagram. 

\medskip

As it is well-known, a {\it framed link} is a pair $(L,c)$, where $L$ is a link in $\mathbb S^3$ with $l \ge 1$ components and $c=(c_1,c_2, \dots, c_l)$, is an $l$-tuple of integers. $(L,c)$ represents - in dimension 3 - the 3-manifold $M^3(L,c)$ obtained from $\mathbb S^3$ by Dehn surgery along $(L,c)$, as well as - in dimension 4 - the (simply-connected) PL 4-manifold $M^4(L,c)$, whose boundary coincides with $M^3(L,c)$, obtained  from $\mathbb D^4$ by adding 2-handles along $(L,c)$. 

Moreover, in virtue of a celebrated result by \cite{Montesinos} and \cite{Laudenbach-Poenaru}, in case $M^3(L,c) = \#_r (\mathbb S^1 \times \mathbb S^2)$ (with $r \ge 0$), then the framed link $(L,c)$ represents also the closed PL 4-manifold  $\overline{ M^4(L,c)}$ obtained from $M^4(L,c)$ by adding  - in a unique way - $r$ 3-handles and a 4-handle.   

However, while it is well-known that every 3-manifold $M^3$ admits a framed link $(L,c)$ so that $M^3=M^3(L,c)$, it is an open question whether or not each closed simply-connected PL 4-manifold $M^4$ may be represented by a suitable framed link (or, 
even more, if $M^4$ admits a so called {\it special} handle decomposition, i.e. a handle decomposition lacking in 1-handles and 3-handles: see \cite[Problem 4.18]{Kirby_1995}, \cite{[M]}, \cite{Casali-Cristofori_RACSAM 2021}). 

As far as general compact PL 4-manifolds (with empty or connected boundary) are concerned, it is necessary to extend the notion of framed link, so to comprehend also 
the case of trivial (i.e. unknotted and unlinked)  {\it dotted components}, which represent 1-handles of the associated handle decomposition of the manifold: 
in this way, any framed link $(L^{(m)},d)$ admitting $m \ge 1$ trivial dotted components - which is properly called a {\it Kirby diagram} - uniquely represents the compact PL 4-manifold $M^4(L^{(m)} ,d)$  
obtained  from $\mathbb D^4$ by adding 1-handles according to the $m$ trivial dotted components and 2-handles along the framed components. 
Note that  the boundary of $M^4(L^{(m)} ,d)$ coincides with $M^3(L,c)$, $(L,c)$ being the framed link obtained from the Kirby diagram $(L^{(m)},d)$ 
by substituting each dotted component with a 0-framed one;  hence, in case $M^3(L,c) = \#_r (\mathbb S^1 \times \mathbb S^2)$ (with $r \ge 0$), the Kirby diagram $(L^{(m)},d)$ 
uniquely represents also the closed  PL 4-manifold $\overline{M^4(L^{(m)} ,d)}$ obtained from $M^4(L^{(m)} ,d)$ 
by adding  - in a unique way - $r$ 3-handles and a 4-handle.   

\medskip
  
In this paper we describe how to obtain algorithmically a regular 5-colored graph representing $M^4(L,c)$ (resp. representing $M^4(L^{(m)} ,d))$ directly ``drawn over" a planar diagram of $(L,c)$ (resp. of  
$(L^{(m)} ,d))$: see Procedure B and Theorem \ref{M4(L,c)} in Section \ref{sec.framed_links} (see Procedure C and Theorem \ref{thm.main(dotted)} in Section \ref{sec.dotted_links}). 
Hence, the algorithms allow to construct explicitly triangulations of the compact 4-manifolds associated to 
framed links and Kirby diagrams\footnote{Indeed, both procedures are going to be implemented in a C++ program, connected to the topological software package Regina (\cite{Burton}): see R.A.~Burke, {\em Triangulating Exotic 4-Manifolds}, in preparation.}.

\medskip

As a consequence, the procedures yield upper bounds for both the invariants {\it  regular genus} and {\it gem-complexity}  of the represented $4$-manifolds.  

As regards framed links, the upper bounds - which significantly improve the ones obtained in  \cite{Casali JKTR2000} - are summarized in the following statement where  
$m_\alpha$ denotes the number  of $\alpha$-colored regions in a chess-board coloration of $L$, by colors $\alpha$ and $\beta$ say, with the convention that the infinite region is $\alpha$-colored;
furthermore, if $w_i$  and $c_i$ denote respectively the writhe and the framing of the $i$-th component of $L,$ (for each $i\in\{1,\ldots,l\}$, $l$ being the number of components of $L$), we set:  
$$\bar t_i = \begin{cases}
                                                       \vert w_i-c_i\vert\quad  \text{if} \ w_i \ne c_i  \\  
                                                      2  \quad \text{otherwise}
                                                      \end{cases}$$

\begin{thm}\label{regular-genus&gem-complexity} 
Let $(L,c)$ be a framed link with $l$ components and $s$ crossings. Then, the following estimation of the regular genus of $M^4(L,c)$ holds:
$$\mathcal G(M^4(L,c))\leq m_{\alpha} + l$$
Moreover, if $L$ is not the trivial knot, then the gem-complexity of $M^4(L,c)$ satisfies the following inequality:
$$k(M^4(L,c))\leq 4s - l + 2\sum_{i=1}^l \bar t_i$$
\end{thm}

As regards Kirby diagrams $(L^{(m)} ,d)$, the estimation for the gem-complexity involves the quantity  $\bar t_i$, defined exactly as in the case of framed links, 
but only for the framed components, while the estimation for the regular genus involves a quantity depending on the construction (i.e. the quantity $u$ appearing in Theorem \ref{thm.main(dotted)}), which can be increased by the number of undercrossings of the  framed components.\footnote{Note that previous work \cite{Casali Compl2004} didn't yield upper bounds for gem-complexity or regular genus, since the combinatorial properties of the obtained 5-colored graph with boundary representing $M^4(L^{(m)} ,d)$ could not be  
``a priori" determined.}  

\begin{thm}\label{regular-genus&gem-complexity(dotted)} 
Let $(L^{(m)},d)$ be a Kirby diagram with $s$ crossings, $l$ components, whose first $m \ge 1$ are dotted,  and $\bar s$  
undercrossings of the  framed components; then, 
$$\mathcal G(M^4(L^{(m)},d))\leq  \  s + \bar s  + (l-m) +1$$ 
Moreover, if $L$ is different from the trivial knot,  
$$k(M^4(L^{(m)},d))\leq \ 2s + 2 \bar s  + 2 m -1 + 2\sum_{i=m+1}^l \bar t_i  $$
\end{thm}

\medskip

Various examples are presented, including infinite families of framed links where the above upper bound for the regular genus turns out to be sharp (Example 1   
and Example 2 in Section \ref{sec.framed_links}). 

Moreover, the process is applied in order to obtain a pair of 5-colored graphs representing an exotic pair of compact PL 4-manifolds (i.e. a pair of  4-manifolds which are TOP-homeomorphic but not PL-homeomorphic), thus opening the search for possibile graph-defined PL invariants distinguishing them  (Example \ref{ex.exotic} in Section \ref{sec.framed_links}, with related Figures \ref{fig.W1} and \ref{fig.W2}).

\medskip  

Note that, although for better understanding the procedure regarding framed links is presented in a separate section of the paper, it is nothing but a particular  
case of the one regarding Kirby diagrams with $m \ge 1 $ dotted components. Hence, if we denote by  $(L^{(*)},d)$ an arbitrary Kirby diagram (possibly without dotted components), we can concisely    
state that the paper shows how to obtain a $5$-colored graph representing the compact $4$-manifold $M^4(L^{(*)},d)$, directly ``drawn over" the Kirby diagram $(L^{(*)} ,d)$. 

Finally, we point out that, if  the associated $3$-manifold is the $3$-sphere, then the obtained $5$-colored graph actually represents the closed $4$-manifold $\overline{M^4(L^{(*)},d)}$, too. Hence, the procedure yields triangulations for any closed (simply-connected) PL 4-manifold admitting handle decompositions without 3-handles.

In the general case of Kirby diagrams representing a closed 4-manifold $\overline{M^4(L^{(*)},d)}$ (i.e., according to \cite{Montesinos}, in case of {\it Heegaard diagrams} for closed $4$-manifolds), we hope soon to be able to extend the above procedure, in order to construct algorithmically - at least in a wide set of situations, when the boundary 3-manifold may be combinatorially recognized as  $\#_r (\mathbb S^1 \times \mathbb S^2)$ (with $r \ge 1$) - a $5$-colored graph representing  $\overline{M^4(L^{(*)},d)}$. 

 \bigskip

\section{Colored graphs representing PL manifolds} \label{prelim}

In this section we will briefly recall some basic notions about the representation of  compact   
PL manifolds by regular colored graphs ({\it crystallization theory}).
For more details we refer to the survey papers \cite{Ferri-Gagliardi-Grasselli} and \cite{Casali-Cristofori-Gagliardi Complutense 2015}. 

From now on, unless otherwise stated, all spaces and maps will be considered in the PL category, and all manifolds will be assumed to be connected and orientable\footnote{Actually all concepts and results exist also, with suitable adaptations, for non-orientable manifolds; however, since the present paper focuses on the relationship between Kirby diagrams and colored graphs,
we will restrict to the orientable case.}. 

Crystallization theory was first developed for closed manifolds; the extension to the case of non-empty boundary, that is more recent, is performed by making use of the wider 
class of singular manifolds.

\begin{defn}\label{singular manifold}  {\em A {\it singular (PL) $n$-manifold} is a closed connected $n$-dimensional polyhedron admitting a simplicial triangulation where the links of vertices 
are closed connected $(n-1)$-manifolds, while  the links of all $h$-simplices of the triangulation with $h > 0$ are (PL) $(n-h-1)$-spheres. Vertices whose links are not PL $(n-1)$-spheres are called {\it singular}.}
\end{defn}

\begin{rem} \label{correspondence-sing-boundary}{\em If $N$ is a singular $n$-manifold, then a compact $n$-manifold $\check N$ is easily obtained by deleting small open neighbourhoods of its singular vertices.
Obviously $N=\check N$ iff $N$ is a closed manifold, otherwise $\check N$ has non-empty boundary (without spherical components).
Conversely, given a compact $n$-manifold $M$, a singular $n$-manifold $\widehat M$ can be constructed by capping off each component of $\partial M$ by a cone over it.

Note that, by restricting ourselves to the class of compact $n$-manifolds with no spherical boundary components,  the above correspondence is bijective and so singular $n$-manifolds and compact $n$-manifolds of this class can be associated  to each other in a well-defined way. 

For this reason, throughout the present work, we will make a further restriction considering only compact manifolds without spherical boundary components. 
Obviously, in this wider context, closed $n$-manifolds are characterized by $M= \widehat M.$}
\end{rem}

\begin{defn} \label{$n+1$-colored graph}
 {\em An {\it $(n+1)$-colored graph}  ($n \ge 2$) is a pair $(\G,\g)$, where $\G=(V(\G), E(\G))$ is a multigraph (i.e. multiple edges are allowed, but no loops) which is regular of degree  $n+1$  (i.e. each vertex has exactly $n+1$ incident edges), and  $\g: E(\G) \rightarrow \Delta_n=\{0,\ldots, n\}$ is a map which is injective on adjacent edges ({\it edge-coloration}).}
\end{defn}

In the following, for sake of concision, when the coloration is clearly understood, we will drop it in the notation for a colored graph.   As usual, we will call {\it order} of a colored graph the number of its vertices.  

\smallskip

For every  $\{c_1, \dots, c_h\} \subseteq\Delta_n$ let $\G_{\{c_1, \dots, c_h\}}$  be the subgraph obtained from $\G$ by deleting all the edges that are not 
colored by $\{c_1, \dots, c_h\}$. 
Furthermore, the complementary set of $\{c\}$ (resp. $\{c_1,\dots,c_h\}$)  in $\Delta_n$ will be denoted by $\hat c$ (resp. $\hat c_1\cdots\hat c_h$). 
The connected components of $\G_{\{c_1, \dots, c_h\}}$ are called {\it $\{c_1, \dots, c_h\}$-residues} of $\G$; 
their number will be denoted by $g_{\{c_1, \dots, c_h\}}$ (or, for short, by $g_{c_1,c_2}$ and $g_{\hat c}$ if $h=2$ and $h = n$ respectively). 
\medskip 

For any $(n+1)$-colored graph $\G$, an $n$-dimensional simplicial cell-complex $K(\G)$ can be constructed in the following way: 
\begin{itemize}
\item the $n$-simplexes of $K(\G)$ are in bijective correspondence with the vertices of $\G$ and each $n$-simplex has its vertices injectively labeled by the elements of $\Delta_n$;
\item if two vertices of $\G$ are $c$-adjacent ($c\in\Delta_n$), the $(n-1)$-dimensional faces of their corresponding $n$-simplices that are opposite to the $c$-labeled vertices are identified, so that equally labeled vertices coincide.
\end{itemize}

\smallskip

$\vert K(\Gamma)\vert$ turns out to be an $n$-pseudomanifold and $\Gamma$ is said to {\it represent} it. 

Note that, by construction, $\G$ can be seen as the 1-skeleton of the dual complex of $K(\G)$.
As a consequence there is a bijection between the $\{c_1, \dots, c_h\}$-residues of  $\G$ and the $(n-h)$-simplices of $K(\G)$ whose vertices are labeled by $\hat c_1 \cdots \hat c_h$. 
In particular, given an $(n+1)$-colored graph $\G$, each connected component of $\G_{\hat c}$ ($c\in\Delta_n$) is an $n$-colored graph representing the disjoint link\footnote{Given a simplicial cell-complex $K$ and an $h$-simplex $\sigma^h$ of $K$, the {\it disjoint star} of $\sigma^h$ in $K$ is the simplicial cell-complex obtained by taking all $n$-simplices of $K$ having $\sigma^h$ as a face and identifying only their faces that do not contain $\sigma^h.$ The {\it disjoint link}, $lkd(\sigma^h,K)$, of $\sigma^h$ in $K$ is the subcomplex of the disjoint star formed by those simplices that do not intersect $\sigma^h.$} of a $c$-labeled vertex of $K(\G)$, that is also (PL) homeomorphic to the link of this vertex in the first barycentric subdivision of $K.$

Therefore, we can characterize $(n+1)$-colored graphs representing singular (resp. closed) $n$-manifolds as satisfying the condition that for each color $c\in\Delta_n$ any $\hat c$-residue represents a connected closed $(n-1)$-manifold\footnote{In case of polyhedra arising from colored graphs, the condition about links of vertices obviously implies the one about links of $h$-simplices, with $h\ge 0.$}  
 (resp. the $(n-1)$-sphere). 
 
\medskip

Furthermore, in virtue of the bijection described in Remark \ref{correspondence-sing-boundary}, an $(n+1)$-colored graph $\G$ is said to {\it represent} a compact $n$-manifold $M$ with no spherical boundary components (or, equivalently, to be a {\it gem} of $M$, where gem means {\it Graph Encoding Manifold})  if $\G$ represents its associated singular manifold, i.e. if $\vert K(\G)\vert=\widehat M$. Actually, if $\partial M\neq\emptyset$, $K(\Gamma)$ naturally gives rise to a  ``triangulation'' of $M$ consisting of partially truncated $n$-simplexes obtained by removing small open neighbourhoods of the singular vertices of $\widehat M$. Therefore, in the present paper,  by a little abuse of notation, we will call $K(\Gamma)$ a {\it triangulation} of $M$ also in the case of a compact manifold with non-empty boundary.

\medskip
The following theorem extends to the boundary case a well-known result - originally stated in \cite{Pezzana}  - founding the combinatorial representation theory for closed manifolds of arbitrary dimension via regular  
colored graphs. 

\begin{thm}{\em (\cite{Casali-Cristofori-Grasselli})}\ \label{Theorem_gem}  
Any compact orientable (resp. non-orientable) $n$-manifold with no spherical boundary components admits a bipartite (resp. non-bipartite) $(n+1)$-colored graph representing it.
\end{thm}

If $\G$ is a gem of a compact $n$-manifold,  an $n$-residue of $\G$ will be called  {\it singular}  if it does not represent $\mathbb S^{n-1}.$ 
Similarly, a color $c$ will be called {\it singular} if at least one of the $\hat c$-residues of $\G$ is singular.

An advantage of colored graphs as representing tools for compact $n$-manifolds is the possibility of combinatorially defining PL invariants.

One of the most important and studied among them is the {\it (generalized) regular genus} extending to higher dimension the classical genus of a surface and the Heegaard genus of a $3$-manifold. 
Spheres are characterized by having null regular genus (\cite{Ferri-Gagliardi Proc AMS 1982}), while classification results according to regular genus and concerning $4$- and $5$-manifolds can be found in \cite{Casali-Gagliardi ProcAMS}, \cite{Casali_Forum2003},  \cite{Casali-Cristofori-Gagliardi Complutense 2015}, \cite{Casali-Cristofori_generalized-genus} both for the closed and for the non-empty boundary case.

The definition of the invariant relies on the following result about the existence of a particular type of embedding of colored graphs into closed surfaces.

\begin{prop}{\em (\cite{Gagliardi 1981})}\label{reg_emb}
Let $\G$ be a bipartite\footnote{Since this paper concerns only orientable manifolds, we have restricted the statement only to the bipartite case, although a similar result holds also for non-bipartite graphs.} $(n+1)$-colored graph of order $2p$. Then for each cyclic permutation $\varepsilon = (\varepsilon_0,\ldots,\varepsilon_n)$ of $\Delta_n$, up to inverse, there exists a cellular embedding, called \emph{regular}, of $\G$  into an orientable closed surface $F_{\varepsilon}(\G)$ whose regions are bounded by the images of the $\{\varepsilon_j,\varepsilon_{j+1}\}$-colored cycles, for each $j \in \mathbb Z_{n+1}$.
Moreover, the genus $\rho_{\varepsilon} (\G)$ of $F_{\varepsilon}(\G)$ satisfies

\begin{equation}\label{reg-genus}
2 - 2\rho_\varepsilon(\G)= \sum_{j\in \mathbb{Z}_{n+1}} g_{\varepsilon_j, \varepsilon_{j+1}} + (1-n)p.
\end{equation}

\end{prop}

\begin{defn}  {\em The {\it regular genus} of an $(n+1)$-colored graph $\G$ is defined as}
$$\rho(\G) = min\{\rho_\varepsilon(\G)\ \vert\ \varepsilon\ \text{cyclic permutation of \ } \Delta_n\}; $$
{\em the  {\it (generalized) regular genus} of a compact $n$-manifold $M$ is defined as}
$$\mathcal G (M) = \min\{\rho(\G)\ \vert\ \G\ \text{gem of \ } M\}.$$
\end{defn}

Within crystallization theory a notion of ``complexity'' of a compact $n$-manifold arises naturally and, similarly to other concepts of complexity (for example Matveev's 
complexity for 3-manifolds) is related to the minimum number of $n$-simplexes in a colored triangulation of the associated singular manifold:

\begin{defn}  {\em The {\it (generalized) gem-complexity} of a compact $n$-manifold $M$ is defined as }
$$k(M) = \min\{p-1\ \vert\ \exists \textit{\ a gem of\ } M\textit{\ with \ }2p \textit{\ vertices}\}$$
\end{defn}

Important tools in crystallization theory are combinatorial moves transforming colored graphs without affecting the represented manifolds (see for example \cite{Ferri-Gagliardi-Grasselli}, \cite{Ferri-Gagliardi Pacific}, \cite{Bandieri-Gagliardi}, \cite{Lins-book}, \cite{Lins-Mulazzani}); we will recall here only the most important ones, while other moves will be introduced in the following sections.  

\begin{defn}  {\em An {\it $r$-dipole ($1\le r\le n$) of colors $c_1,\ldots,c_r$} in an $(n+1)$-colored graph $\G$ is a subgraph of $\G$ consisting in two vertices joined by $r$ edges, colored by $c_1,\ldots,c_r$, such that the vertices belong to different $\hat c_1\ldots\hat c_r$-residues of $\G$. 
An $r$-dipole can be eliminated from $\G$ by deleting the subgraph and welding the remaining hanging edges according to their colors; in this way another $(n+1)$-colored graph $\G^\prime$ is obtained. The addition of the dipole to $\G^\prime$ is the inverse operation.\\
The dipole is called {\it proper} if $\vert K(\G)\vert$ and $\vert K(\G^\prime)\vert$ are 
 (PL) homeomorphic. } \end{defn}
    
\begin{prop}\label{proper-dipole}
{\rm (\cite[Proposition 5.3]{Gagliardi 1987})} An $r$-dipole ($1\le r\le n$) of colors $c_1,\ldots,c_r$ in an $(n+1)$-colored graph $\G$ is proper if and only if at least one of the two connected components of $\G_{\hat c_1\ldots\hat c_r}$ intersecting the dipole represents the $(n-r)$-sphere.
\end{prop}

\noindent Without going into details, we point out that - as proved in the quoted paper - the elimination (or the addition)  of a proper dipole corresponds to a re-triangulation of a suitable ball embedded in the cell-complex associated to the colored graph.  

\begin{rem} {\em Note that, if $\G$ represents a compact $n$-manifold $M$, then all $r$-dipoles with $1 < r \le n$ are proper; 
further, if $M$ has either empty or connected boundary, then $1$-dipoles are proper, too.}
\end{rem}

Given an arbitrary $(n+1)$-colored graph representing a compact $n$-manifold $M$ with empty or connected boundary, then by eliminating all possible (proper) $1$-dipoles, we can always obtain an $(n+1)$-colored graph $\G$ still representing $M$ and such that for each color $c\in\Delta_n$, $\G_{\hat c}$ is connected. 
Such a colored graph is called a {\it crystallization} of $M.$
Moreover, it is always possible to assume - up to permutation of the color set - that any gem (and, in particular, any crystallization) of such a manifold, has color $n$ as its (unique) possible singular color.  

\smallskip

Finally,  as already hinted to in the Introduction, we recall that a  graph-based representation for compact PL manifolds with non-empty boundary - different from the one considered in this section - was already introduced by Gagliardi in the eighties (see \cite{Ferri-Gagliardi-Grasselli}) by means of colored graphs failing to be regular.   

More precisely, any compact  $n$-manifold can be represented by a pair $(\Lambda, \lambda)$, where $\lambda$ is still an edge-coloration on $E(\Lambda)$ by means of $\Delta_n$, but $\Lambda$ may miss some (or even all) $n$-colored edges: such a $(\Lambda, \lambda)$ is said to be an {\it $(n+1)$-colored graph with boundary, 
regular with respect to color $n$}, and vertices missing the $n$-colored edge are called {\it boundary vertices}.

However, a connection between these different kinds of representation can be established through an easy combinatorial procedure, called {\it capping-off}. 

\begin{prop} {\em (\cite{Ferri-Gagliardi Yokohama 1985})} 
\label{cappingoff} Let $(\Lambda, \lambda)$ be an $(n+1)$-colored graph with boundary, regular with respect to color $n$, representing the compact $n$-manifold $M$.  
Chosen a color $c \in \Delta_{n-1}$, let $(\G,\g)$ be the regular $(n+1)$-colored graph obtained from $\Lambda$ by \emph{capping-off with respect to color $c$}, i.e. by joining  two boundary vertices by an $n$-colored edge, whenever they belong to the same $\{c,n\}$-colored path in $\Lambda.$ Then,  $(\G,\gamma)$ represents 
the singular $n$-manifold $\widehat M$, and hence $M$, too.
\end{prop}

\section{From framed links to $5$-colored graphs}  \label{sec.framed_links}

In this section we will present a construction that enables to obtain $5$-colored graphs representing all compact (simply-connected) $4$-manifolds associated to framed links, i.e. Kirby diagrams without dotted components. Note that such a class of compact $4$-manifolds contains also each closed  (simply-connected) $4$-manifold admitting a {\it special} handle decomposition (\cite[Section 3.3]{[M]}), i.e. a handle decomposition containing no $1$- and $3$-handles. 
\medskip

As already recalled in the Introduction, for each framed link $(L,c)$ ($c=(c_1, \dots, c_l)$, with $c_i \in \mathbb Z$ $\forall i \in \{1,\dots, l\}$, $l$ being the number of components of $L$), we denote by $M^4(L,c)$ 
the $4$-manifold with boundary obtained from $\mathbb D^4$ by adding $l$ 2-handles according to the framed link  $(L,c)$.  
The boundary of $M^4(L,c)$ is the closed orientable 3-manifold $M^3(L,c)$ obtained from $\mathbb S^3$ by Dehn surgery along $(L,c)$.  
In case $M^3(L,c)\cong \mathbb S^3$, we will consider, and still denote by $M^4(L,c),$ the closed $4$-manifold obtained by adding a further $4$-handle. 

Now, let us suppose that the link $L$ is embedded in $S^3 = \mathbb R^3\cup\{\infty\}$ so that it admits a regular projection $\pi\ : \mathbb S^3\to\mathbb R^2\times\{0\};$ 
in the following we will identify $L$ with its planar diagram $\pi(L)$, thus referring to {\it arcs, crossings} and {\it regions} of $L$ instead of $\pi(L).$   

Similarly, by the {\it writhe} of a component $L_i$ of $L$ (denoted by $w(L_i)$) we mean the writhe of the corresponding component of $\pi(L).$
For each $i\in\{1,\ldots,l\}$, we say that $L_i$ needs $\vert c_i-w(L_i)\vert$ ``additional curls", which are positive or negative according to whether $c_i$ is
greater or less than $w(L_i)$ (see Figure \ref{curls}).

In \cite{Casali JKTR2000} a construction is described, yielding a $4$-colored graph representing the $3$-manifold associated to a given framed link . 
The procedure consists of the following steps. 

\bigskip

\noindent
{\bf  PROCEDURE A -  from $\mathbf{(L,c)}$ to $\mathbf{\Lambda(L,c)}$  representing $\mathbf{M^3(L,c)}$:}

\begin{enumerate}
\item Each crossing of $L$ gives rise to the order eight graph in Figure \ref{crossing}, while each possible 
 (whether already in $L$ or additional)  curl gives rise to one of the order four  graphs of Figure \ref{graph-curls}-left or Figure \ref{graph-curls}-right according to the curl being positive or negative.
 
 \begin{figure}[!ht] 
\centerline{\scalebox{0.35}{\includegraphics{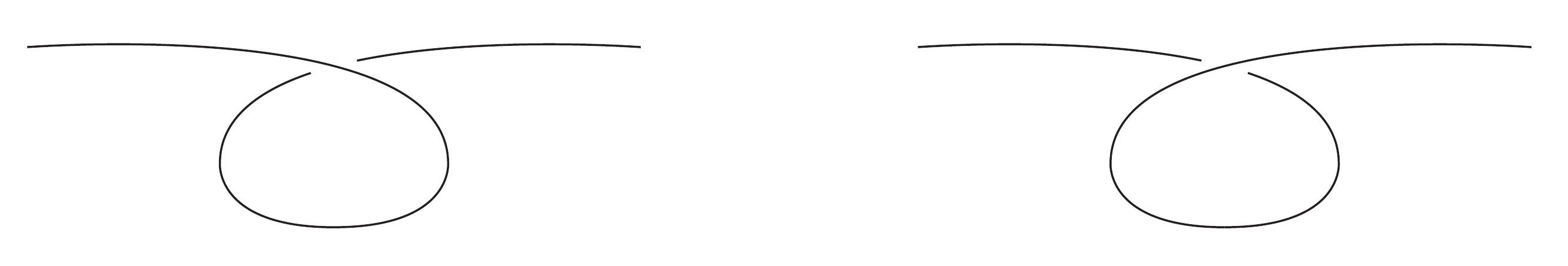}}}
\caption{{\footnotesize Positive (left) and negative (right) curls}}
\label{curls}
\end{figure}
 
 \begin{figure}[!ht] 
\centerline{\scalebox{0.65}{\includegraphics{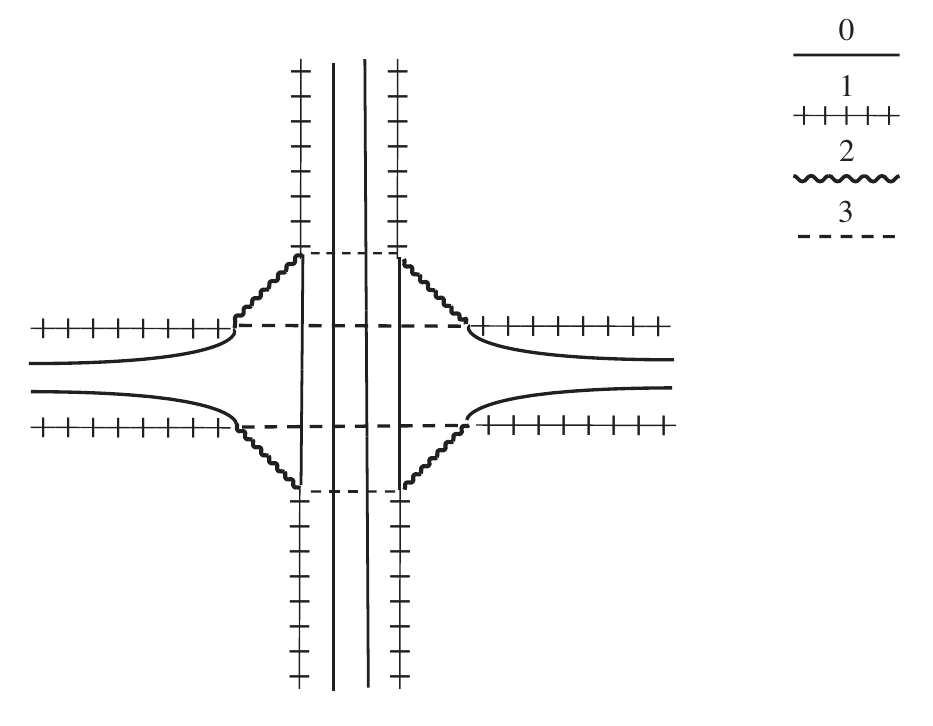}}}
\caption{{\footnotesize  4-colored graph corresponding to a crossing}}
\label{crossing}
\end{figure}

\begin{figure}[!ht] 
\centerline{\scalebox{0.50}{\includegraphics{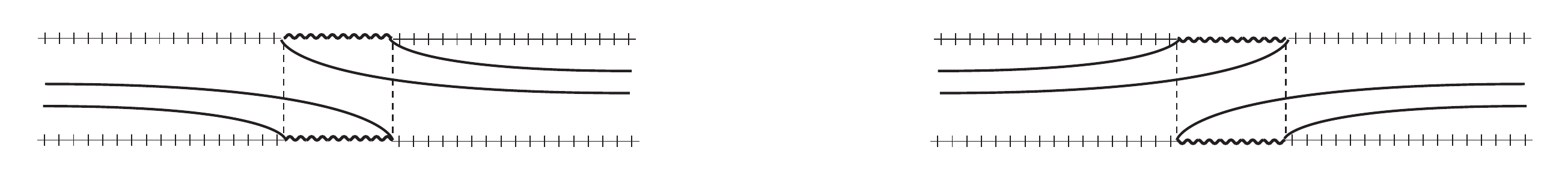}}}
\caption{{\footnotesize 4-colored graphs corresponding to a positive curl (left) and a negative curl (right)}}
\label{graph-curls}
\end{figure}

\item 
The hanging $0$- and $1$-colored edges of the above graphs should be ``pasted"  together so that every region of $L$, having $r$ crossings on its boundary, gives rise to a $\{1,2\}$-colored cycle of length $2r$  (with each 1-colored edge corresponding to a part of the boundary between two crossings) while each component  $L_i$ ($i\in\{1,\ldots,l\}$), having $s_i$ crossings and $t_i$ additional curls, gives rise to two $\{0,3\}$-colored cycles of length $2(s_i+t_i).$
\end{enumerate}

\begin{figure}[!ht] 
\centerline{\scalebox{0.6}{\includegraphics{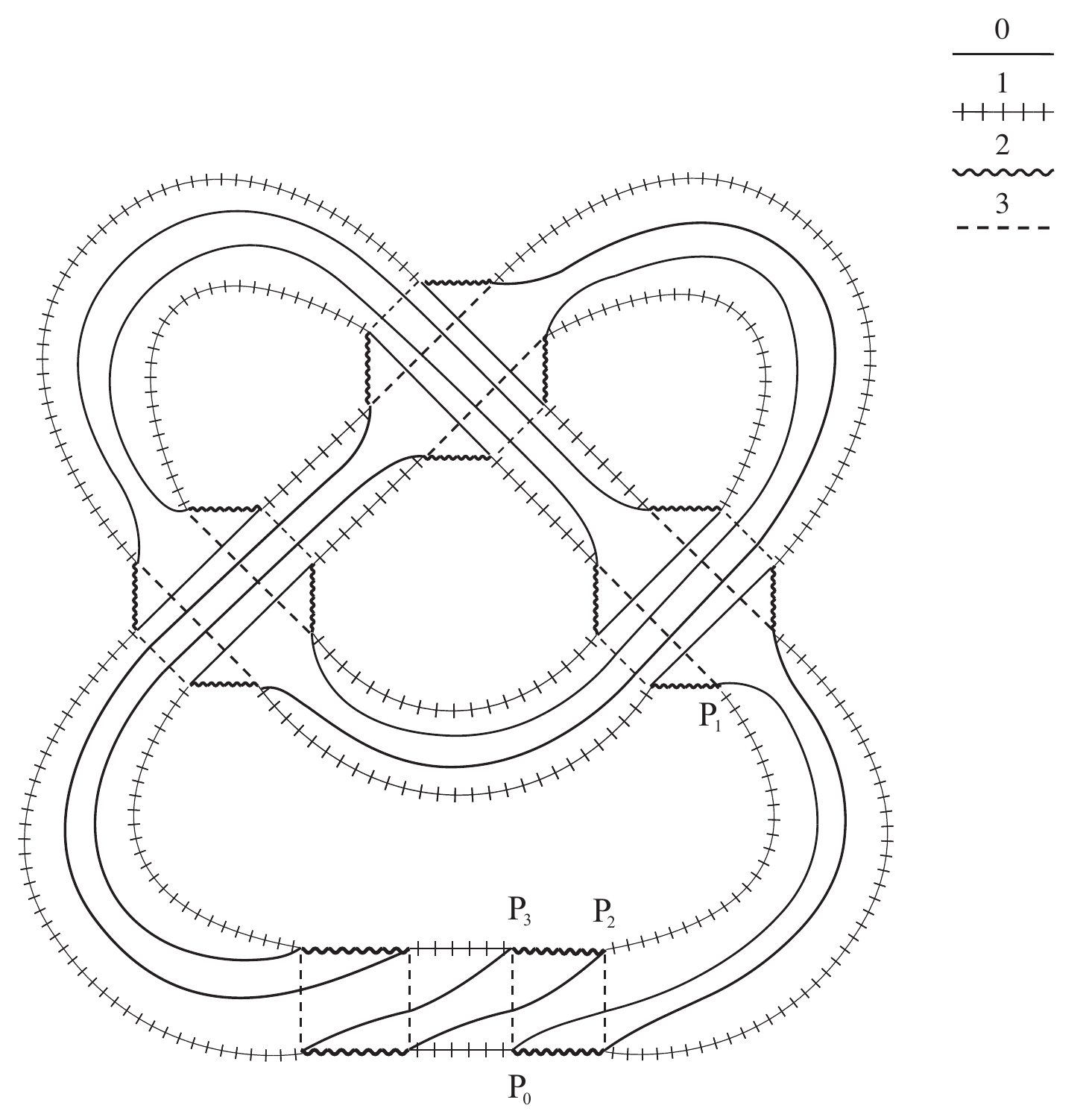}}}
\caption{{\footnotesize The $4$-colored graph $\Lambda(L,c)$ representing  $M^3(L,c)$, for $c=+1$ and $L=$ trefoil}}
\label{fig.trefoil1}
\end{figure}

\begin{rem} \label{rem.drawn_over} 
 {\em As pointed out in \cite{Casali JKTR2000}, $\Lambda(L,c)$ can be directly  ``drawn over'' $L$ (see for example Figure \ref{fig.trefoil1}, obtained by applying Procedure A to  the trefoil knot, with framing $c=+1$). 
In particular, if $a$ is the part of an arc of $L$ lying between two adjacent crossings, there are exactly two $1$-colored edges of $\Lambda(L,c)$ that are ``parallel" to $a$, one for each region of $L$ having $a$ on its boundary.

\noindent 
Moreover, note that - by possibly adding to $L$ a trivial pair of opposite additional curls - a particular subgraph $Q_i$, called {\it quadricolor}, can be selected in $\Lambda(L,c)$  for each component $L_i$ of $L$ ($i \in \{1, \dots, l\}$).
A quadricolor consists of four vertices $\{P_0,\ P_1,\ P_2,\ P_3\}$ such that $P_s,P_{s+1}$ are connected by an $s$-colored edge (for each $s\in\mathbb Z_4$)  and $P_s$ does not belong to the $\{s+1,s+2\}$-colored cycle shared by the other three vertices. 
It is not difficult to see that - in virtue of the above described procedure A - such a situation arises with $\{P_0,\ P_2,\ P_3\}$ belonging to the subgraph corresponding to a curl and $P_1$ to an adjacent undercrossing or curl of the same sign (see again Figure \ref{fig.trefoil1}, where the vertices of the quadricolor are highlighted).}
\end{rem}

Let us now describe how to construct, starting from a given framed link, a 5-colored graph  which will be proved to represent the $4$-manifold associated to the framed link itself. 

\bigskip

\noindent 
{\bf  PROCEDURE B -  from $\mathbf{(L,c)}$ to $\mathbf{\Gamma(L,c)}$ \ (representing $\mathbf{M^4(L,c)}$):}

\begin{figure}
\centerline{\scalebox{0.6}{\includegraphics{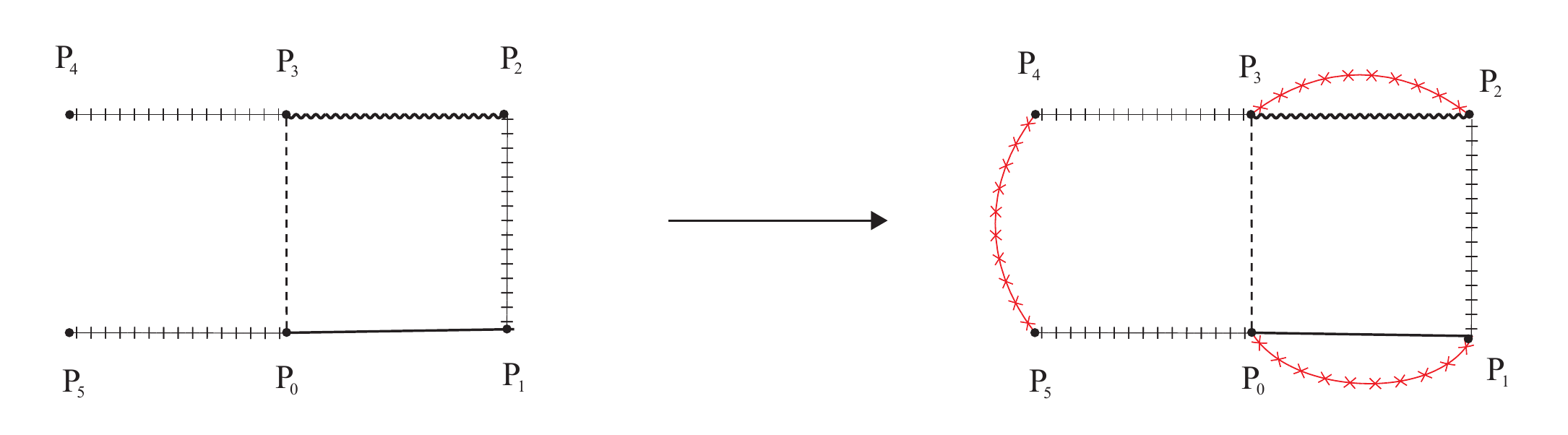}}}
\caption{{\footnotesize  Main step yielding $\G(L,c)$}}
\label{fig.quadricolor-singular}
\end{figure}

\begin{enumerate}
\item   Let  $\Lambda(L,c)$  be the 4-colored graph constructed from $(L,c)$ according to Procedure $A$.
\item  For each component $L_i$ of $L$ ($i \in \{1, \dots, l\}$), choose a quadricolor $Q_i$, according to Remark  \ref{rem.drawn_over}. 
For each $i \in \{1, \dots, l\}$, add a triad of $4$-colored edges between the vertices $P_{2r}$ and $P_{2r+1}$, $\forall r \in \{0,1,2\}$, involved in the quadricolor $\mathcal Q_i$ (as shown in Figure \ref{fig.quadricolor-singular}). 
\item Add $4$-colored edges between the remaining vertices of $\Lambda(L,c)$, so to ``double" the $1$-colored ones.
\end{enumerate}
\bigskip

\begin{figure}
\centerline{\scalebox{0.6}{\includegraphics{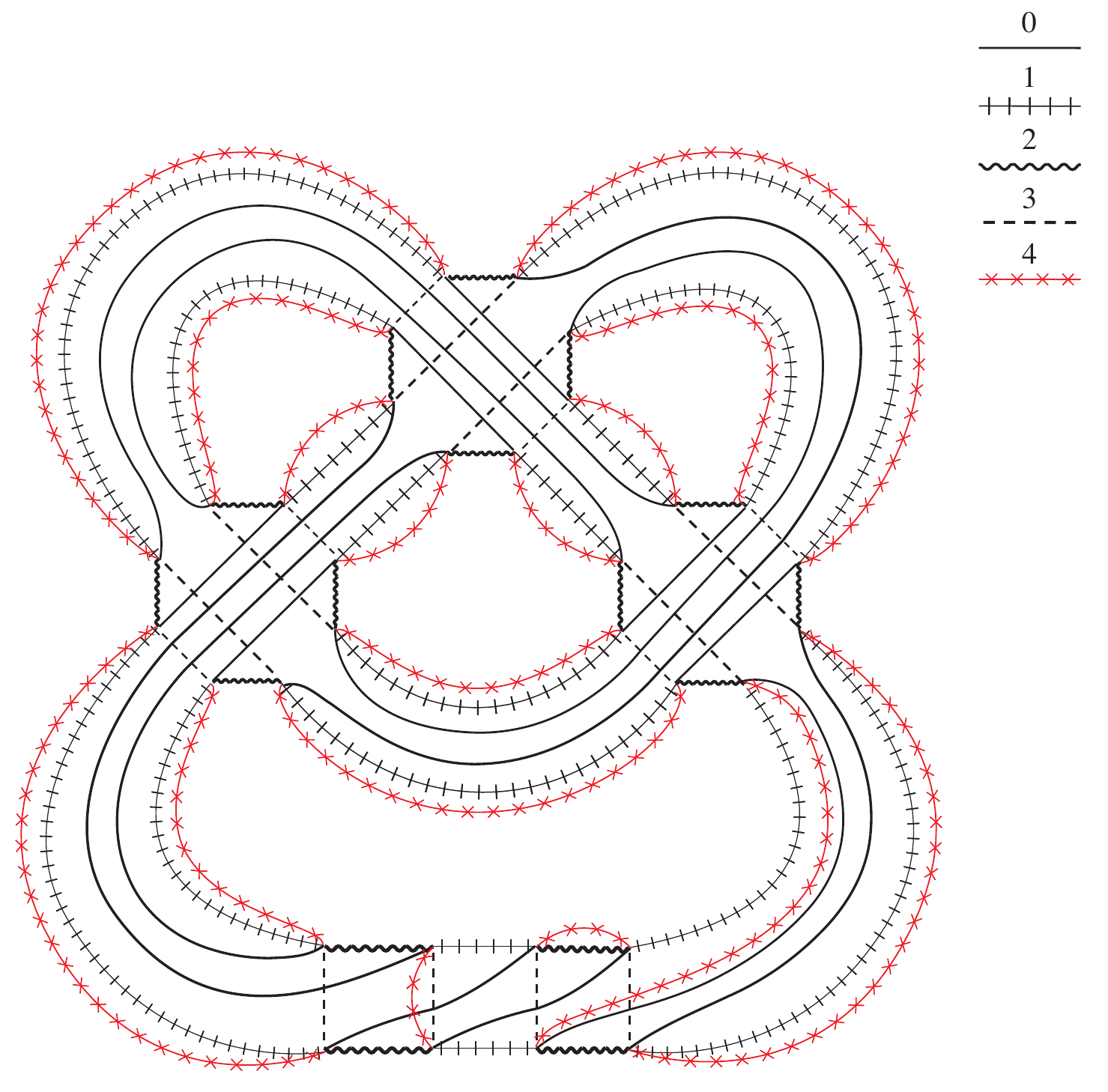}}}
\caption{{\footnotesize The $5$-colored graph $\Gamma(L,c)$, for $c=+1$ and $L=$ trefoil} }
\label{fig.trefoil2}
\end{figure}

Figure \ref{fig.trefoil2} shows the $5$-colored graph $\G(L,c)$ in the case of the trefoil knot with framing $+1$.

\bigskip

The following theorem states  that - as already disclosed  
- the $5$-colored graph $\Gamma(L,c)$ represents $M^4(L,c)$.
Moreover, the theorem also states the existence of a further 5-colored graph representing $M^4(L,c)$, with reduced regular genus, whose estimation involves the number $m_\alpha$ 
of $\alpha$-colored regions in a chess-board coloration of $L$, by colors $\alpha$ and $\beta$ say, with the convention that the infinite region is $\alpha$-colored. 

With this aim, if $(L,c)$ is a framed link with $l$ components and $s$ crossings, let us recall that, 
for each $i\in\{1,\ldots,l\}$, we set $\bar t_i = \begin{cases}
                                                       \vert w_i-c_i\vert\quad  if\ w_i \ne c_i  \\  
                                                      2  \quad otherwise
                                                      \end{cases}$ where $w_i$ denotes the writhe of the $i$-th component of $L.$ 

\medskip

\begin{thm}\label{M4(L,c)} \ 
\begin{itemize}
 \item [(i)] 
For each framed link $(L,c)$, the $5$-colored graph $\G(L,c)$ obtained via Procedure B represents the compact $4$-manifold $M^4(L,c);$ it has regular genus less or equal to $s+l+1$ and, if $L$ is different from the trivial knot,\footnote{More precisely we suppose the projection $\pi(L)$ to be different from the standard diagram of the trivial knot. This case, which is already well-known (see \cite{Casali-Cristofori_generalized-genus}), is nevertheless discussed in details in Example \ref{example1}.} its order is \ $8s + 4\sum_{i=1}^l \bar t_i;$
 \item [(ii)]  via a standard sequence of graph moves, a $5$-colored graph, still representing  $M^4(L,c)$, can be obtained, 
 whose regular genus is less or equal to $m_{\alpha} + l$, while the regular genus of its $\hat 4$-residue, representing $\partial M^4(L,c)=M^3(L,c),$ is less or equal to $m_\alpha.$
\end{itemize}

\end{thm}

Theorem \ref{M4(L,c)} will be proved in Section \ref{sec.proofs}.

\medskip

As a direct consequence of Theorem \ref{M4(L,c)}, upper bounds can be established for both the invariants regular genus and gem-complexity of a compact $4$-manifold 
represented by a framed link $(L,c)$, in terms of the combinatorial properties of the link itself, as already stated in Theorem \ref{regular-genus&gem-complexity} in the Introduction.

\medskip

\noindent\textit{Proof of Theorem~{\upshape\ref{regular-genus&gem-complexity}}} 
The upper bound for the regular genus of  $M^4(L,c)$  trivially follows from Theorem {\ref{M4(L,c)}} (ii).  

As regards the upper bound for the gem-complexity of  $M^4(L,c)$, we have to make use of the computation of the order of  $\G(L,c)$ obtained in Theorem {\ref{M4(L,c)}} (i), 
but also to note that - as already pointed out in \cite{Casali JKTR2000} - the 4-colored graph $\Lambda(L,c)$ has exactly  
$l$ $\hat 2$-residues, and that the same happens for $\G(L,c)$;  hence, by deleting $l-1$     
(proper) $1$-dipoles, a new $5$-colored graph $\G^{\prime}(L,c)$ representing  $M^4(L,c)$ may be obtained, with 
$$ \#V(\G^{\prime}(L,c)) =  \#V(\G(L,c)) - 2(l-1)  =  8s - 2l + 2 + 4\sum_{i=1}^l \bar t_i  $$\qed

The case of the trivial knot is discussed in the following example.

\begin{example}\label{example1} {\em Let $(K_0,c)$ be the trivial knot with framing $c\in\mathbb N\cup\{0\}$; 
if $c\geq 2$, then $K_0$ requires $c$ additional positive curls and the $5$-colored graph $\Gamma(K_0,c)$, with $4c$ vertices, which is obtained by applying Procedure B, turns out to coincide with the one that in \cite{Casali-Cristofori_generalized-genus} is 
proved to represent exactly $\xi_c$, the $\mathbb D^2$-bundle over $\mathbb S^2$ with Euler number $c$,
as expected from Theorem \ref{M4(L,c)}. Furthermore, if $c$ is even, it is known that $k(L(c,1))=2c-1$ (see \cite[Remark 4.5]{Casali-Cristofori lens spaces}); hence 
$\Gamma(K_0,c)$ realizes the gem-complexity of $\xi_c$, and therefore the second bound of Theorem \ref{M4(L,c)} is sharp.
 
If $c=0$ (resp. $c=1$), then $K_0$ requires two positive and two negative (resp. two positive and one negative) additional curls in order to get a quadricolor;
however in this case the resulting graph $\Gamma(K_0,c)$ admits a sequence of dipole moves consisting in three $3$-dipoles and one $2$-dipole 
(resp. consisting in two $3$-dipoles) cancellations yielding a minimal order eight crystallization of $\mathbb S^2\times\mathbb D^2$ 
(resp. the minimal order eight crystallization of $\mathbb{CP}^2$) obtained in \cite{Casali-Cristofori_generalized-genus} (resp. in \cite{Gagliardi CP2}). 
 
Note that for each $c\in\mathbb N\cup\{0\}$, $\Gamma(K_0,c)$ realizes the regular genus of the represented 4-manifold, which is equal to $2$ ($=m_{\alpha} +l$), as proved in \cite{Casali-Cristofori_generalized-genus}.  
Hence, for this infinite family of compact 4-manifolds, the first upper bound of Theorem \ref{regular-genus&gem-complexity} turns out to be sharp.}
 \end{example}

We will end this section with further examples of the described construction.

\begin{example} {\em Let $(L_H,c)$ be the Hopf link and $c=(\bar c,0)$ with $\bar c$ even (resp. odd); then Procedure B yields a $5$-colored graph that, by Theorem \ref{M4(L,c)}, represents $\mathbb S^2\times\mathbb S^2$ (resp. $\mathbb{CP}^2\#(-\mathbb{CP}^2)$) and realizes its regular genus (which is known to be equal to 4: see \cite{Casali-Cristofori-Gagliardi Complutense 2015} and references therein). 
In particular, if $\bar c=0$, a sequence of dipole cancellations and a {\it $\rho$-pair switching} (see Definition \ref{rho-pair} in Section \ref{sec.dotted_links}), applied to $\Gamma(L_H,(0,0))$, yield a $5$-colored graph which belongs to the existing catalogue\footnote{More details about such catalogue (together with other similar ones) can be found at\\ 
https://cdm.unimore.it/home/matematica/casali.mariarita/CATALOGUES.htm\#dimension\_4} of crystallizations of $4$-manifolds up to gem-complexity 8 (see \cite{Casali-Cristofori ElecJComb 2015}).}
\end{example}

\begin{example} {\em Let $M^4(L,c)$ be a linear plumbing of spheres, whose boundary is therefore the lens space $L(p,q)$ such that $(c_1,\ldots,c_l)$ is a continued fraction expansion
of $-\frac pq$ (\cite{GS}); then, by Theorem \ref{regular-genus&gem-complexity}, the regular genus of $M^4(L,c)$ is less or equal to $2l.$  }
\end{example}

\begin{figure}[!h]
\centerline{\scalebox{0.45}{\includegraphics{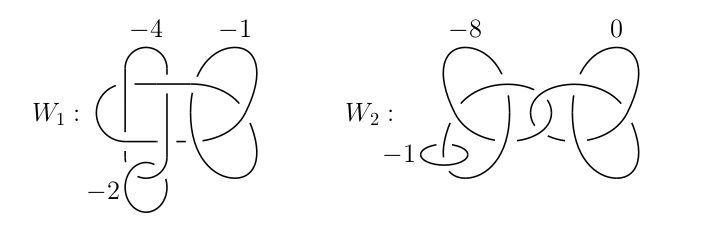}}}
\caption{{\footnotesize Framed links representing the exotic pair $W_1$ and $W_2$ (pictures from \cite{Naoe_ExperMath})}}
\label{exotic_pair}
\end{figure}

\begin{figure}[!h]
\centerline{\scalebox{0.31}  {\includegraphics{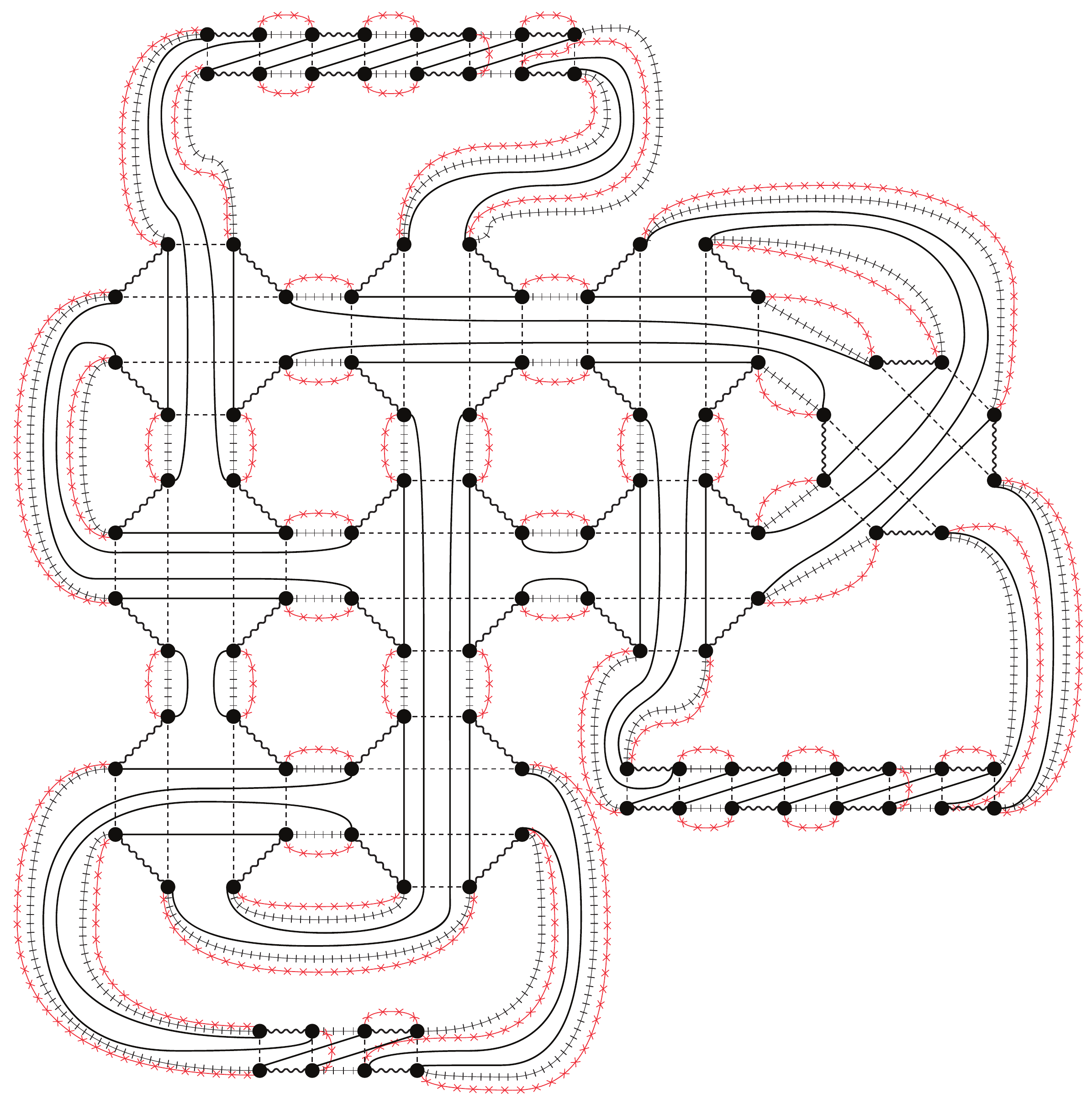}}} 
\caption{{\footnotesize A 5-colored graph representing the compact simply-connected PL 4-manifold $W_1$}}
\label{fig.W1}
\end{figure}

\begin{figure}[!h]
\centerline{\scalebox{0.36}    {\includegraphics{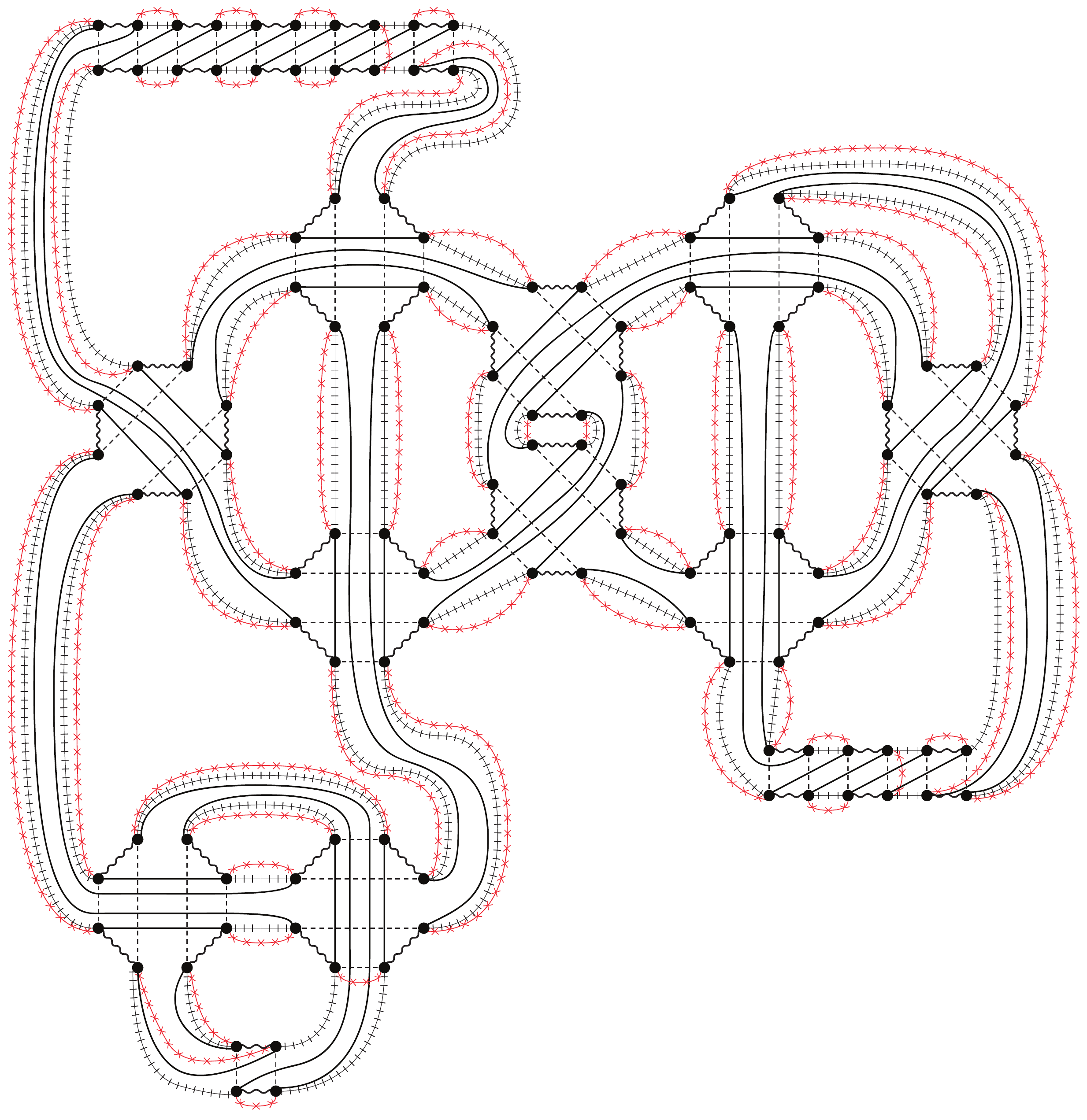}}} 
\caption{{\footnotesize A 5-colored graph representing the compact simply-connected PL 4-manifold $W_2$}}
\label{fig.W2}
\end{figure}

\begin{example} \label{ex.exotic} {\em Procedure B, applied to the framed links description given in \cite{Naoe_ExperMath} of an exotic pair (see Figure \ref{exotic_pair}), allows to obtain two regular 5-colored graphs representing two compact simply-connected PL 4-manifolds $W_1$ and $W_2$ with the same topological structure that are not PL-homeomorphic: see Figures \ref{fig.W1} and \ref{fig.W2}, which obviously encode two triangulations of $W_1$ and $W_2$ respectively.

Other applications of the procedures obtained in the present paper, in order to get triangulations of exotic 4-manifolds, will appear in R.A.~Burke, {\em Triangulating Exotic 4-Manifolds} (in preparation).}
\end{example}

\section{Proof of Theorem \ref{M4(L,c)}}\label{sec.proofs}    

Roughly speaking, the proof of  the first statement of Theorem \ref{M4(L,c)} - i.e. the fact 
that $\Gamma(L,c)$ represents $M^4(L,c)$ -  will be performed by means of the followings steps: 
\begin{itemize}
\item[(i)] starting from the $4$-colored graph $\Lambda(L,c)$ - already proved to represent $M^3(L,c)$ in \cite{Casali JKTR2000} - we obtain a $4$-colored graph $\Lambda_{smooth}$ representing $\mathbb S^3$ by suitably exchanging a triad of $1$-colored edges for each component of $L$ (Proposition \ref{Gamma_smooth}(i)); 
\item[(ii)] by capping-off with respect to color 1, we obtain a 5-colored graph representing  $\mathbb D^4$;  
\item[(iii)] by re-establishing the triads of $1$-colored edges, the 5-colored graph $\Gamma(L,c)$  is obtained. Since the only singular 4-residue of $\Gamma(L,c)$ is the $\hat 4$-residue $\Lambda(L,c)$,  $\Gamma(L,c)$  represents a 4-manifold with connected boundary $M^3(L,c)$; moreover, $\Gamma(L,c)$ represents $M^4(L,c)$ since each triad exchanging is proved to correspond to the addition of a 2-handle according to the related framed component (Proposition \ref{Gamma_smooth}(ii)). 
\end{itemize}

\medskip

\noindent Let us now go into details. \\  
Given a framed link $(L,c)$, we can always assume that, for each component $L_i$ ($i\in\{1,\ldots,l\}$), an additional curl is placed near an undercrossing;
as observed in Section \ref{sec.framed_links} such a configuration gives rise, in the $4$-colored graph $\Lambda(L,c)$, to a quadricolor that we denote by $\mathcal Q_i.$

By cancelling the quadricolor $\mathcal Q_i$ and pasting the resulting hanging edges of the same color, we obtain a new $4$-colored graph $\Lambda^{(\hat\imath)}(L,c)$;
we call this operation the \textit{smoothing} of the quadricolor $\mathcal Q_i.$

The following proposition shows that the smoothing of a quadricolor in a $4$-colored graph obtained from a framed link via Procedure B (see Section 
\ref{sec.framed_links}) turns out to be equivalent to the Dehn surgery on the complementary knot of the involved link component.
More precisely, with the above notations, the result can be stated as follows:

\begin{prop}\label{smoothing-quadricolor} If $\Lambda^{(\hat\imath)}(L,c)$ is the 4-colored graph obtained from $\Lambda(L,c)$ by smoothing the quadricolor of 
the $i$-th framed component, then $K(\Lambda^{(\hat\imath)}(L,c))$ is obtained from $K(\Lambda(L,c))$ by Dehn surgery on the complementary knot of $L_i.$ 

Hence, $K(\Lambda^{(\hat\imath)}(L,c))$ represents the 3-manifold associated to the framed link $(L^{\hat\imath}, c^{\hat\imath})$ obtained from $(L, c)$ by deleting the $i$-th 
component.
 \end{prop}

\dimo Let $(L^{(\tilde\imath)},c^{(\tilde\imath)})$ denote the $l+1$ components link obtained from $L$ by adding the {\it complementary knot of $L_i$}, i.e. a framed 0 trivial knot linking the component $L_i$ geometrically once; 
moreover, let us suppose that the added trivial component is inserted between the curl and the crossing corresponding to the quadricolor $\mathcal Q_i.$ 
Then, let us consider the $4$-colored graph $\Lambda(L^{(\tilde\imath)},c^{(\tilde\imath)})$ obtained by applying Procedure B  
of Section \ref{sec.framed_links} to the framed link $(L^{(\tilde\imath)},c^{(\tilde\imath)}).$

\begin{figure}[!h]
\centerline{\scalebox{0.6}{\includegraphics{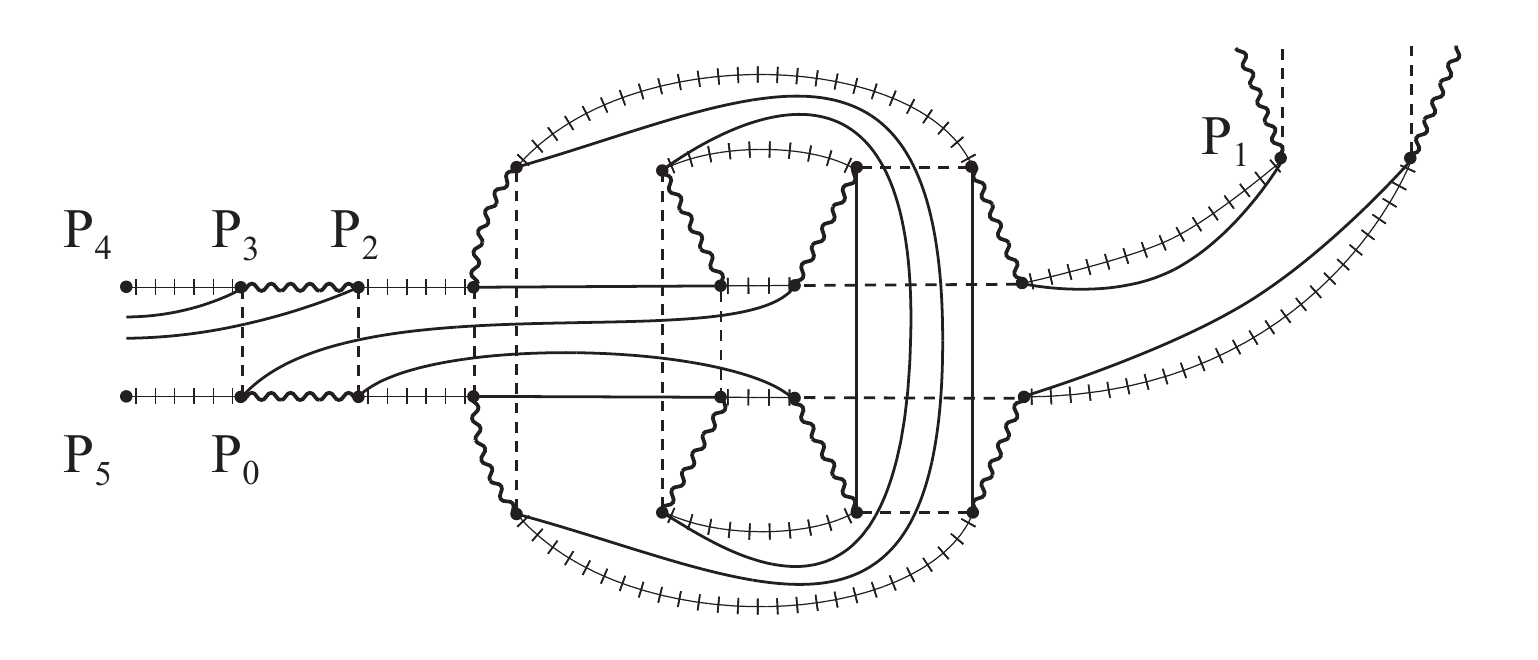}}}
\caption{}
\label{fig.lemmaJKTRinizio}
\end{figure}

\begin{figure}[!h]
\centerline{\scalebox{0.6}{\includegraphics{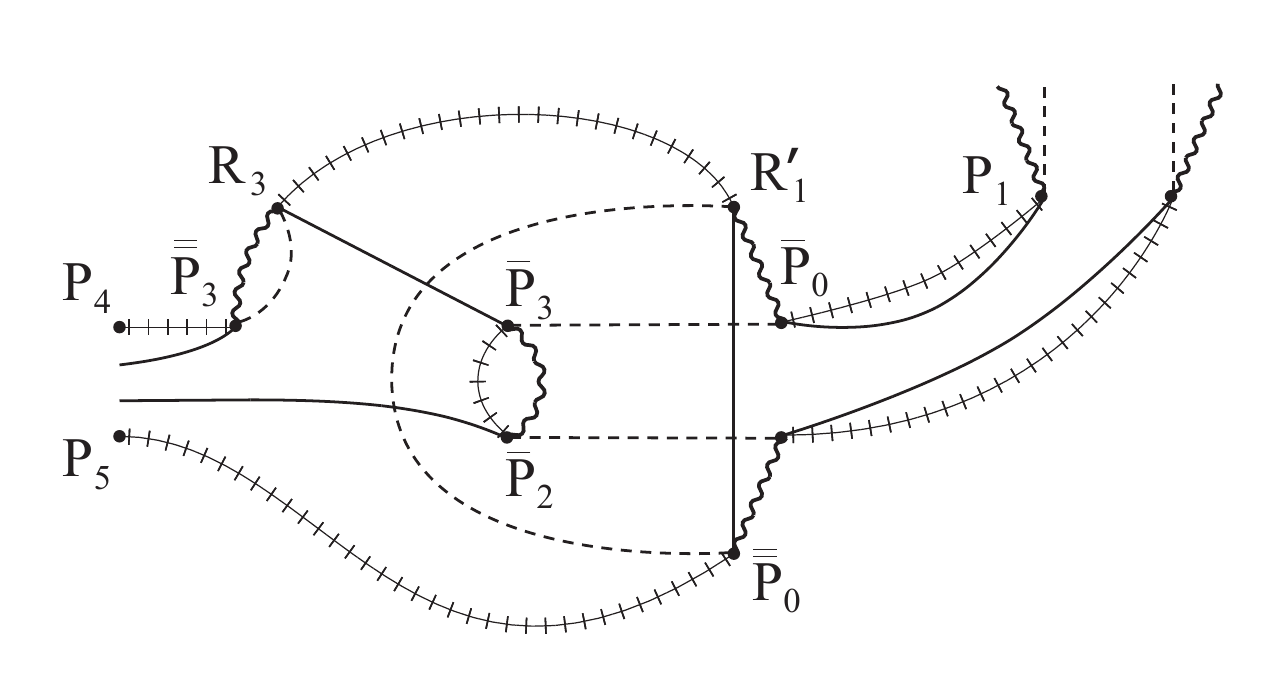}}} 
\caption{}
\label{fig.lemmaJKTRfine}
\end{figure}

$\Lambda(L^{(\tilde\imath)},c^{(\tilde\imath)})$  is everywhere like $\Lambda(L,c)$ except ``near" the quadricolor $\mathcal Q_i$, where it contains the subgraph in Figure \ref{fig.lemmaJKTRinizio} (we denote by $P_j,\  j\in\{0,1,2,3\}$ the vertices of $\mathcal Q_i$, even if they are no longer a quadricolor in $\Lambda(L^{(\tilde\imath)},c^{(\tilde\imath)})$).
In the proof of Lemma 4 of \cite{Casali JKTR2000} it is shown that the above subgraph yields, through a sequence of eliminations of dipoles, the subgraph in Figure \ref{fig.lemmaJKTRfine}.

By subsequently cancelling the 2-dipoles of vertices $\{\bar {\bar P}_0,R'_1\},$ $\ \{\bar P_0,P_1 \},$ $\ \{\bar P_2, \bar P_3 \},$ $\ \{\bar {\bar P}_3,R_3\}$, all vertices of the quadricolor $\mathcal Q_i$ are eliminated and the obtained $4$-colored graph  is precisely $\Lambda^{(\hat\imath)}(L,c).$

Since the addition to $(L,c)$ of the complementary knot of $L_i$ corresponds to the Dehn surgery along it, the  first part of the statement is proved. 
Moreover, the last part follows directly by noting that the component $L_i$ and its complementary knot constitute a pair of complementary handles, whose cancellation does not affect the represented 3-manifold.\qed

\begin{rem}  {\em Quadricolors in $4$-colored graphs were originally introduced by Lins. 
Note that the transformation from $\Lambda(L,c)$ to the $4$-colored graph where the quadricolor $\mathcal Q_i$ is replaced by the subgraph in Figure \ref{fig.lemmaJKTRfine} corresponds to the substitution, in the pseudocomplex $K(\Lambda(L,c))$, of a solid torus with another solid torus having the same boundary. Hence, as already observed by Lins himself, the smoothing of a quadricolor in any $4$-colored graph is equivalent to perform a Dehn surgery on the represented manifold. The above proposition allows, when considering $4$-colored graphs arising from framed links, to identify this surgery precisely as the one along the complementary knot of the component naturally associated to the quadricolor.}
\end{rem}

\begin{figure}[!h]
\centerline{\scalebox{0.55}{\includegraphics{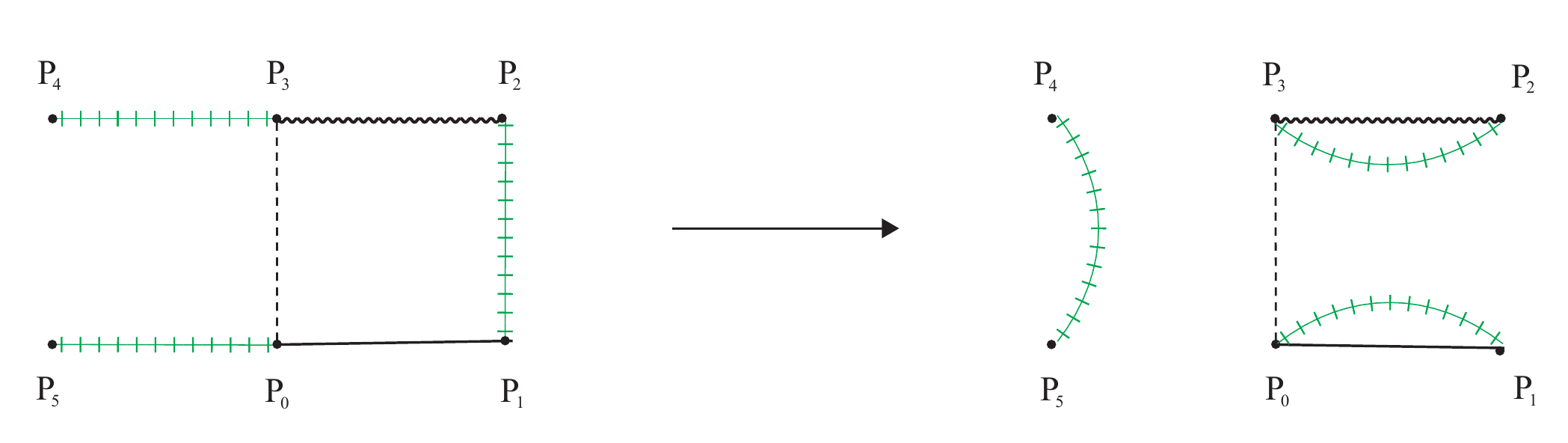}}}
\caption{{\footnotesize Exchanging the triad of 1-colored edges in a quadricolor (I)}}
\label{fig.Gamma_smooth}
\end{figure}

\begin{figure}[!h]
\centerline{\scalebox{0.55}{\includegraphics{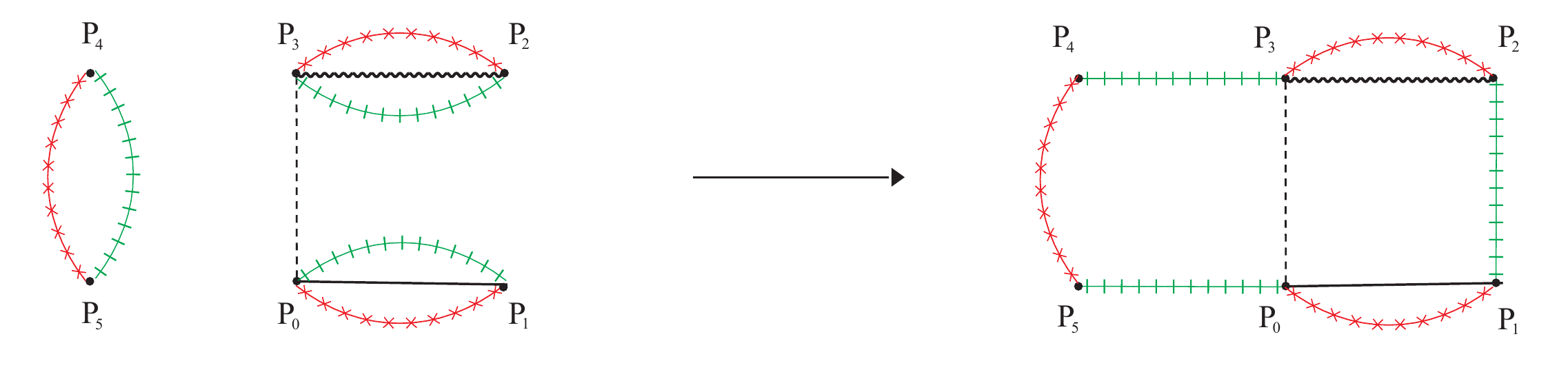}}}
\caption{{\footnotesize Exchanging the triad of 1-colored edges in a quadricolor (II)}}
\label{fig.2-handle}
\end{figure}

\begin{prop}\label{Gamma_smooth} \  \par 
\begin{itemize}
\item[(i)]
The 4-colored graph $\Gamma^{(\hat \imath)}_{smooth}$ (resp. $\Gamma_{smooth}$), obtained from $\Lambda(L,c)$ by 
exchanging the triad of $1$-colored edges (according to Figure \ref{fig.Gamma_smooth}) in a quadricolor of the $i$-th component of $(L,c)$ (resp. in a quadricolor for each framed component of $(L,c)$), represents  the 3-manifold associated to the framed link $(L^{\hat\imath}, c^{\hat\imath})$ obtained from $(L, c)$ by deleting the $i$-th component (resp. represents $\mathbb S^3$).  
\item[(ii)] Let $\tilde\Gamma_{smooth}$ be the 5-colored graph obtained from $\Gamma_{smooth}$ by ``capping off" with respect to color $1$. 
Then, the 5-colored graph $\tilde \Gamma^{(i)}$, obtained from $\tilde \Gamma_{smooth}$ by exchanging the triad of $1$-colored edges (according to Figure \ref{fig.2-handle}) in a quadricolor of the $i$-th component of $(L,c)$, represents  the 4-manifold obtained from $\mathbb D^4$  by adding a 2-handle according to the $i$-the component of $(L,c))$. 
\end{itemize}
 \end{prop}

\dimo 

(i) \ It is sufficient to make use of the proof of Proposition \ref{smoothing-quadricolor} and to note that $\Gamma^{(\hat \imath)}_{smooth}$ is obtained  (modulo the name exchange of $P_j$ into $\bar P_j$, for  $j \in \{0,2,3\}$) from the subgraph in Figure \ref{fig.lemmaJKTRfine} via cancellation of the two 2-dipoles of vertices $\{\bar {\bar P}_0,R'_1\},\  
\{\bar {\bar P}_3,R_3 \}$ in the quadricolor of the $i$-th framed component, while $\Gamma_{smooth}$ is obtained by performing the same procedure for each component of $(L,c)$. 

\begin{figure}
\centerline{\scalebox{0.55}{\includegraphics{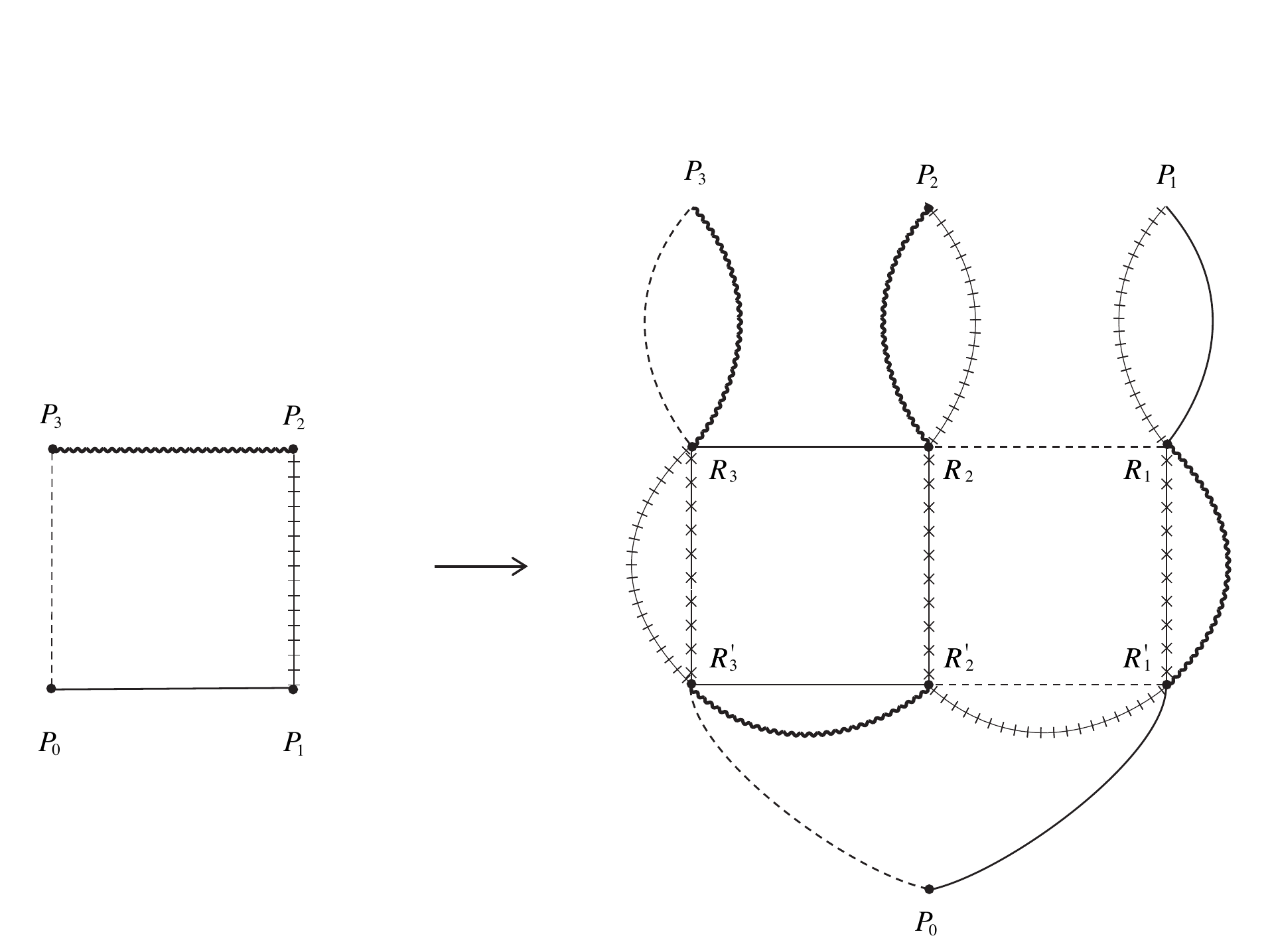}}}
\caption{{\footnotesize Main step from  $\Lambda(L,c)$  to  $\tilde \Lambda(L,c)$  }}
\label{fig.quadricolor}
\end{figure}

(ii) \ It is easy to check that, by a standard sequence of dipole addition, the 4-colored graph $\Gamma_{smooth}$ may be transformed (modulo the name exchange of $P_2$ into $\bar P_2$ and $P_j$ into $\bar{\bar P}_j$, for  $j \in \{0,3\} $) in the 4-colored graph $\tilde \Lambda_{\hat 4}(L,c)$, already considered both in  \cite{Casali JKTR2000} and in  \cite{Casali Compl2004}: more precisely, for each component of the link, it is necessary to add the 2-dipoles of vertices $\{\bar {\bar P}_0,R'_1\}$ and $\{\bar {\bar P}_3,R_3 \}$  shown in Figure \ref{fig.lemmaJKTRfine}, and then to add a 2-dipole of vertices $\{R'_2,R'_3\}$ within the $1$-colored edge with endpoints $\{R'_1,R_3\}$ and the $3$-colored edge with endpoints $\{\bar {\bar P}_0,R'_1\}$. 
The 4-colored graph $\tilde \Lambda_{\hat 4}(L,c)$ is deducible from Figure \ref{fig.quadricolor}, that illustrates the main step  to obtain the 5-colored graph with boundary $\tilde \Lambda(L,c)$, representing $M^4(L,c)$,  from the 4-colored graph  $\Lambda(L,c)$.  

\begin{figure}
\centerline{\scalebox{0.54}{\includegraphics{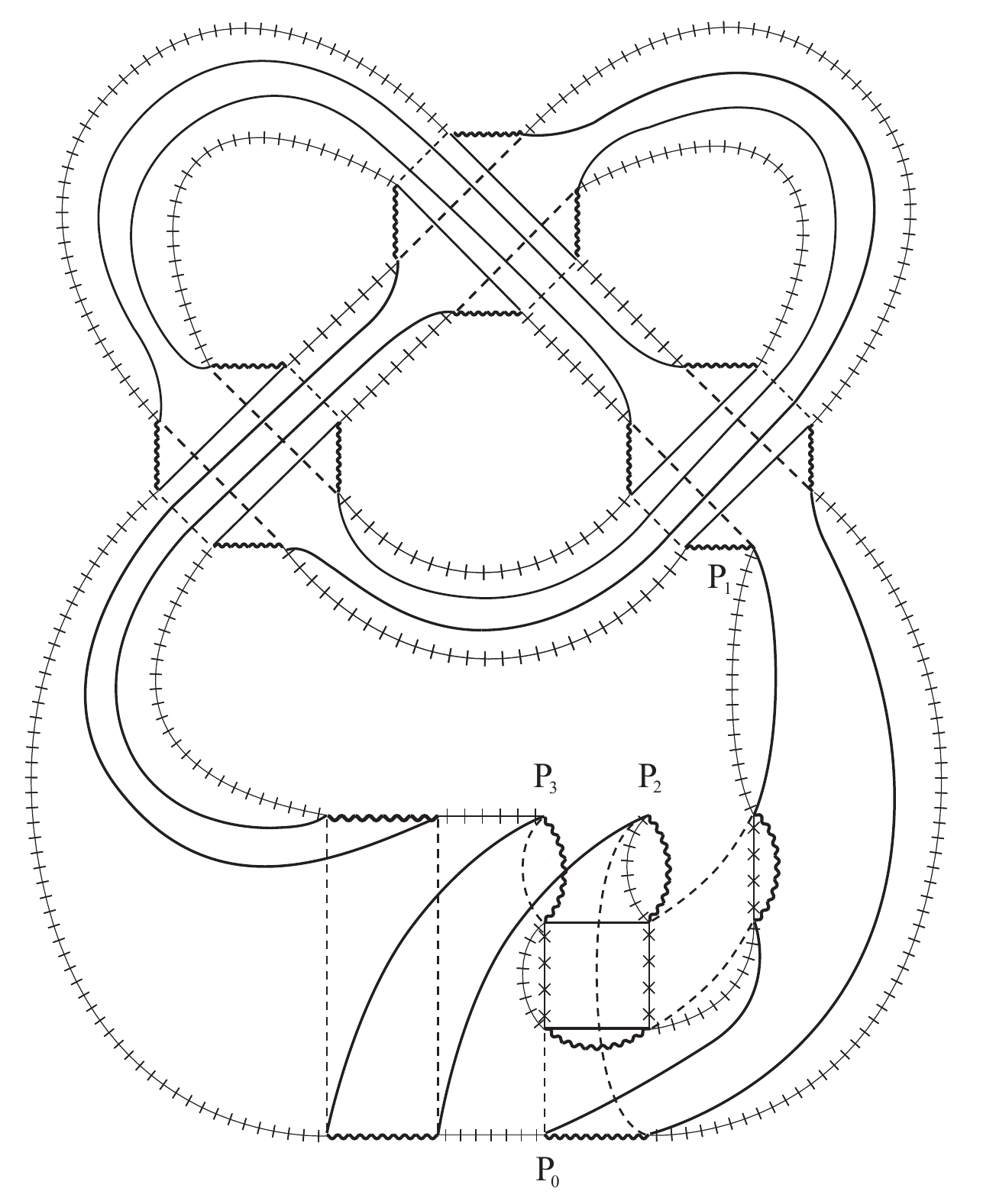}}}
\caption{{\footnotesize The $5$-colored graph with boundary $\tilde \Lambda(L,c)$ representing  $M^4(L,c)$, for $c=+1$ and $L=$ trefoil.}}
\label{fig.M4(trefoil1)}
\end{figure}

Moreover, as explained in \cite[pp. 442-443]{Casali Compl2004},  the 1-skeleton of the associated colored triangulation $K(L,c)= K(\tilde \Lambda_{\hat 4})$ of $\mathbb S^3$, contains two copies $L^{\prime}=L^{\prime}_1 \cup \dots \cup L^{\prime}_l$ and $L^{\prime \prime}=L^{\prime\prime}_1 \cup \dots \cup L^{\prime\prime}_l$  of $L$, with linking number $c_i$ 
between $L^{\prime}_i$ and $L^{\prime\prime}_i$, for each $i \in \{1, \dots, l\}$, and the addition of the triad of $4$-colored edges with endpoints $\{R_1,R_1^{\prime}\}$, $\{R_2,R_2^{\prime}\}$, $\{R_3,R_3^{\prime}\}$ corresponds - as proved in  \cite[Theorem 3]{Casali JKTR2000} - to the attachment on the boundary of $\mathbb D^4$ (i.e. the cone over $K(L,c)$) of a 2-handle whose attaching map sends $L_i^{\prime}$ into $L_i^{\prime\prime}$ (see Figure \ref{fig.quadricolor}-right, and Figure \ref{fig.M4(trefoil1)} for an example of the $5$-colored graph with boundary\footnote{Recall that in this type of colored graphs, some vertices lack incident $4$-colored edges.}  $\tilde \Lambda(L,c)$, where  $(L,c)$ is the trefoil knot with framing $+1$).

Now, if the ``capping off" procedure described  in Proposition \ref{cappingoff} is applied with respect to color $1$, the obtained regular 5-colored graph  (which represents the compact 4-manifold obtained from $\mathbb D^4$ by adding a 2-handle along the $i$-th component of $(L,c)$) turns out to admit a sequence of three 2-dipoles involving only vertices of the quadricolor and never involving color $4$: in fact, they consist of the pairs of vertices $\{P_3,R_3\}$, $\{R_3^{\prime},R_2^{\prime}\}$, $\{P_0,R_1^{\prime}\}$ in Figure \ref{fig.quadricolor} (right). 
It is not difficult to see that, after these cancellations, we obtain exactly the (regular) $5$-colored graph $\tilde \Gamma^{(i)}$, obtained from $\tilde \Gamma_{smooth}$ (which obviously represents $\mathbb D^4$) by cyclically exchanging the triad of $1$-colored edges (according to Figure \ref{fig.2-handle}) in the quadricolor $Q_i$ of the $i$-th component of $(L,c)$. 
 \qed

\begin{rem} \label{rem.sequence_dipoles}   
 {\em Note that a standard sequence of dipole moves exists, transforming $\Gamma^{(\hat \imath)}_{smooth}$ into $\Lambda(L^{\hat\imath}, c^{\hat\imath})$: it follows $L_i$ starting from the quadricolor $Q_i$, where the triad of $1$-colored edges have been cyclically exchanged as in Figure \ref{fig.Gamma_smooth} (right), by deleting first the 2-dipoles $\{P_2, P_3\}$ and $\{P_4, P_5\}$, and then all subsequently generated 2-dipoles, among pairs of vertices, belonging to different bipartition classes, which are either endpoints of $1$-colored edges ``parallel" to adjacent segments of $L_i$, or $0$-adjacent vertices of the subgraph associated to an undercrossing of $L_i$, till to obtain an order two component of the 4-colored graph consisting only of the vertex $P_0$ and its $2$-adjacent vertex.   
Obviously, if the procedure is applied for each $i\in \{1, \dots, l\}$, a standard sequence of dipole eliminations is obtained, transforming the 4-colored graphs  $\Gamma_{smooth}$ into the order two 4-colored graph representing $\mathbb S^3$. }
\end{rem}

\noindent\textit{Proof of Theorem~{\upshape\ref{M4(L,c)}}}

\noindent (i)\ In order to prove that $\Gamma(L,c)$ represents $M^4(L,c)$, it is sufficient to note that the main step yielding $\Gamma(L,c)$, depicted in Figure \ref{fig.quadricolor-singular}, exactly coincides  with the transformation from the 4-colored graph of Figure \ref{fig.Gamma_smooth}-left (representing $M^3(L,c)$) to the 5-colored graph of Figure  
\ref{fig.2-handle}-right (representing $M^4(L,c)$, if the procedure is applied to a quadricolor for each component of $(L,c)$). Hence, Proposition \ref{Gamma_smooth} (i) and (ii) ensure  $\Gamma(L,c)$ actually to represent the compact 4-manifold obtained from $\mathbb D^4$ by adding $l$ 2-handles according to the $l$ components of $(L,c)$. 

\noindent  Now note that, by construction, the $4$-colored graph $\Lambda(L,c)$ has $8s + 4\sum_{i=1}^l \vert w_i-c_i\vert$ vertices. As already observed, the presence of a curl near an undercrossing in a component of $L$ yields a quadricolor $Q_i$.
Therefore, for each $i\in\{1,\ldots,l\}$, if $\vert w_i-c_i\vert\neq 0$, then the required addition of curls ensures the existence of a quadricolor relative to $L_i$, while if $\vert w_i-c_i\vert = 0$ a pair of opposite curls has to be added in order to produce one  (except in the case of the trivial knot which is discussed in Example \ref{example1}).
Since each curl contributes with $4$ vertices to the final $5$-colored graph, the statement concerning the order of $\Gamma(L,c)$  is easily proved. 

With regard to the regular genus of $\Gamma(L,c)$, let us consider the cyclic permutations $\bar\varepsilon=(1,0,2,3)$ of $\Delta_3$ and $\varepsilon=(1,0,2,3,4)$ of $\Delta_4.$

As already pointed out in \cite{Casali JKTR2000}, the construction of $\Lambda(L,c)$ directly yields $\rho_{\bar\varepsilon}(\Lambda(L,c))=s+1.$ On the other hand, it is easy to check - via formula \eqref{reg-genus} - that 
$$2\rho_{\varepsilon}(\Gamma(L,c))- 2\rho_{\bar\varepsilon}(\Lambda(L,c)) = g_{1,3}-g_{3,4}-g_{1,4}+p,$$
\noindent  where $2p$ is the order of $\Gamma(L,c)$ (as well as of $\Lambda(L,c)$).

Since, by construction, $g_{3,4}=g_{1,3}$ and $g_{1,4}=p-2l$, we obtain: 
\begin{equation}\label{difference_reggenus}
\rho_{\varepsilon}(\Gamma(L,c))= \rho_{\bar\varepsilon}(\Lambda(L,c)) +l.
\end{equation}
The result about the regular genus of $\Gamma(L,c)$ now directly follows.

\medskip

\noindent (ii)\ 
As proved in \cite[Theorem 1]{Casali JKTR2000}, the $4$-colored graph $\Lambda(L,c)$ admits a finite sequence of moves, called {\it generalized dipole eliminations}\footnote{A generalized dipole in a $4$-colored graph representing a closed $3$-manifold is a particular subgraph, whose cancellation factorizes into a sequence of proper dipole moves; from the topological point of view, this move corresponds to a Singer move of type III' involving a pair of curves in a suitable Heegaard diagram which can be associated to the $4$-colored graph (see \cite{Ferri-Gagliardi Pacific} for details).}, which preserve the represented manifold and do not affect the quadricolor structures, but reduce the regular genus.
Hence, a new 4-colored graph  $\Omega(L,c)$ representing $M^3(L,c)$ is obtained, having regular genus $m_\alpha$ with respect to the cyclic permutation $\bar\varepsilon = (1,0,2,3)$ of $\Delta_3$ (see  \cite{Casali JKTR2000} for details).  $\Omega(L,c)$ contains a quadricolor for each component of $L$, too, and the results of Proposition  \ref{Gamma_smooth} (i) and (ii) may be applied, exactly as previously done for  $\Lambda(L,c)$, so to obtain - via the move depicted in Figure \ref{fig.quadricolor-singular} performed on a quadricolor for each component of $L$ -  a new 5-colored graph $\tilde \Gamma(L,c)$ representing $M^4(L,c)$.\footnote{It is not difficult to check that $\tilde \Gamma(L,c)$ could also be obtained through the $5$-colored graph with boundary $\tilde \Omega(L,c)$, constructed in \cite{Casali JKTR2000} by applying the move depicted in Figure \ref{fig.quadricolor} for each component of the link: in fact, $\tilde \Omega(L,c)$ represents $M^4(L,c)$, too, and in order to obtain $\tilde \Gamma(L,c)$ it is sufficient to make the capping off  with respect to color $1$ and to delete three $2$-dipoles for each quadricolor, exactly as done in the proof of Proposition \ref{Gamma_smooth} (ii) for $\tilde \Gamma_{smooth}$.}

Now, it is not difficult to check that - in full analogy to equation \eqref{difference_reggenus} - the following relation holds between the regular genera of  $\tilde \Gamma(L,c)$ and $\Omega(L,c)$, with respect to $\bar\varepsilon$ and $\varepsilon$ respectively:
$$ \rho_{\varepsilon}(\tilde\Gamma(L,c))= \rho_{\bar\varepsilon}(\Omega(L,c)) +l.$$
Then, both statements of Theorem \ref{M4(L,c)}(ii) directly follow from $\rho_{\bar\varepsilon}(\Omega(L,c))=m_\alpha$: 
$\rho_{\varepsilon}(\tilde\Gamma(L,c))= m_\alpha + l$, while  $\rho_{\bar\varepsilon}((\tilde\Gamma(L,c))_{\hat 4}))= m_\alpha$ (since the $\hat 4$-residue of $\tilde\Gamma(L,c)$ is exactly $ \Omega(L,c)$). 
\qed

\begin{figure}
\centerline{\scalebox{0.65}{\includegraphics{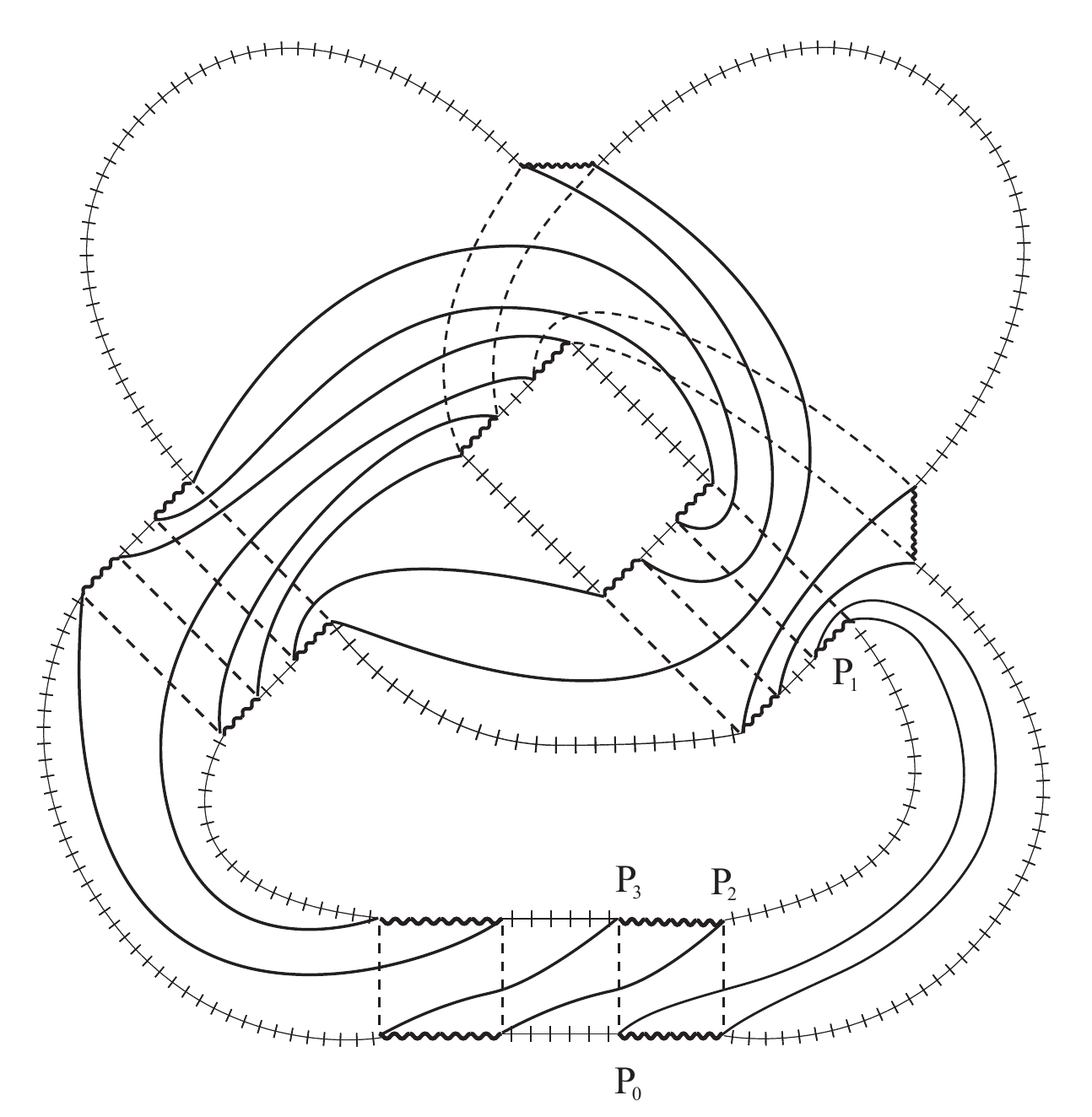}}}
\caption{{\footnotesize The $4$-colored graph $\Omega(L,c)$ representing  $M^3(L,c)$, for $c=+1$ and $L=$ trefoil  }
\label{fig.trefoil1-reduced}  }
\end{figure}

\medskip
See Figures \ref{fig.trefoil1-reduced}  and  \ref{fig.M4(K,c)-reduced} for examples of graphs $\Omega(L,c)$ and $\tilde\G(L,c)$ respectively, where $(L,c)$ is the trefoil knot with framing $+1.$

\begin{figure}
\centerline{\scalebox{0.65}{\includegraphics{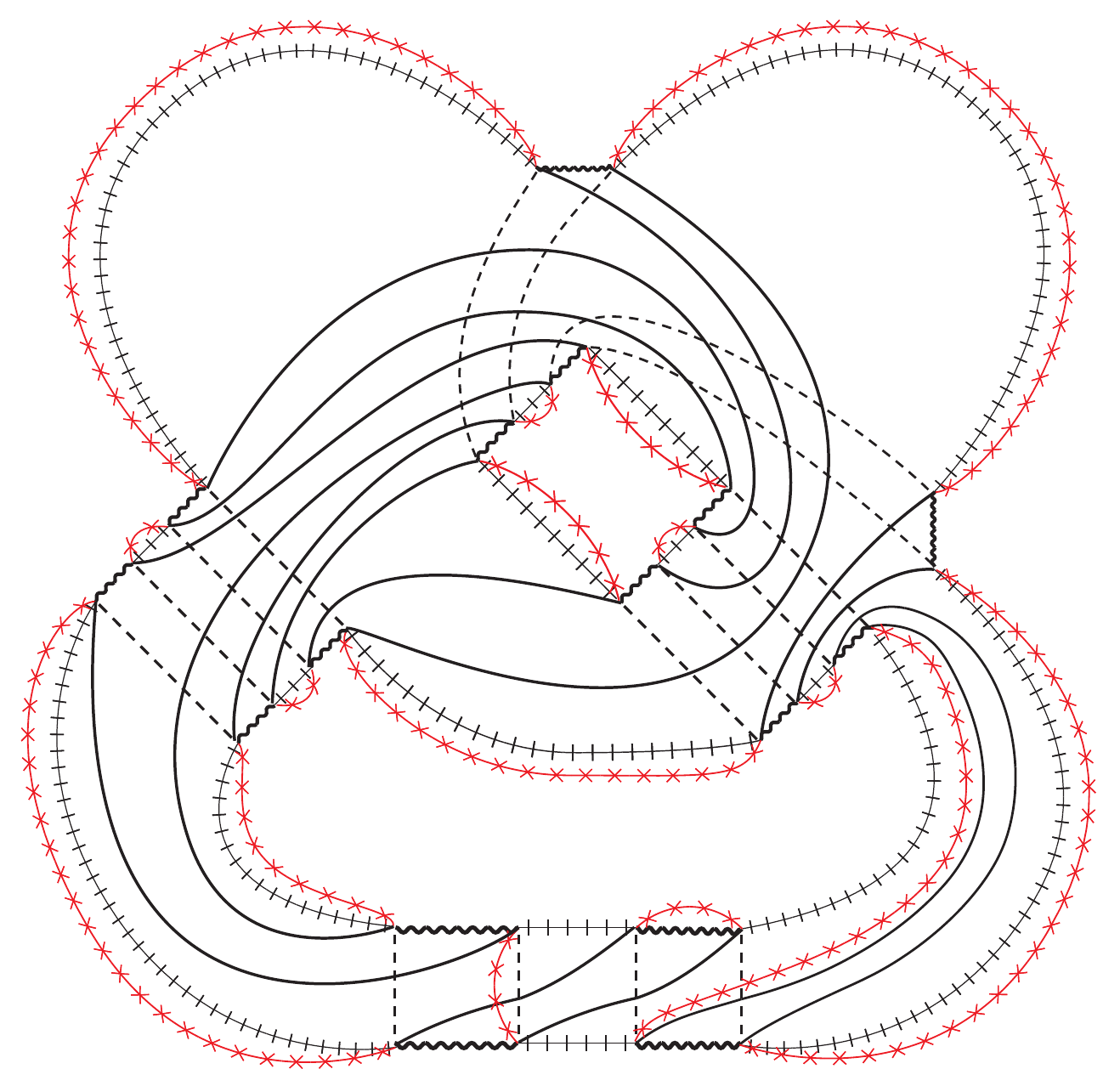}}}
\caption{{\footnotesize The $5$-colored graph $\tilde\G(L,c)$  representing  $M^4(L,c)$, for $c=+1$ and $L=$ trefoil}}
\label{fig.M4(K,c)-reduced}
\end{figure}

\bigskip

\section{From dotted links to $5$-colored graphs  \label{sec.dotted_links}} 

In this section we will take into account the more general case of Kirby diagrams with dotted components, extending the 
procedure and the results of Section  \ref{sec.framed_links}.
Note that, as a consequence, the class of manifolds involved in the construction includes all closed (simply-connected) $4$-manifolds admitting a handle decomposition without 3-handles (\cite[Problem 4.18]{Kirby_1995}, which is of particular interest with regard to exotic PL 4-manifolds: see, for example, \cite{[Ak1]} and \cite{[Ak2]}). 
\medskip

Let $(L^{(m)},d)$ be a Kirby diagram, where $L$ is a link with $l$ components, $L_i$ with $i\in\{1,\ldots,m\}$ (resp. $i\in\{m+1,\ldots,l\}$) being a dotted (resp. framed) component, and $d=(d_1, \dots, d_{l-m})$, where $d_i \in \mathbb Z$ $\forall i \in \{1,\dots, l-m\}$ is the framing of the $(m+i)$-th (framed) component.

As already recalled in the Introduction, we denote by $M^4(L^{(m)},d)$ the $4$-manifold with boundary obtained from $\mathbb D^4$ by adding $m$ 1-handles according to the dotted components and $l-m$ 2-handles according to the framed components of $(L^{(m)},d)$.  
The boundary of $M^4(L^{(m)},d)$ is the closed orientable 3-manifold $M^3(L,c)$ obtained from $\mathbb S^3$ by Dehn surgery along the associated framed link $(L,c)$, obtained by substituting each dotted component by a 0-framed one, i.e. $c=(c_1 \dots, c_{l})$, where $c_i = \begin{cases} 0 & 1\le i \le m \\
                                                                          d_{i-m} & m+1 \le i \le l \end{cases}$

In case $M^3(L,c)\cong \mathbb S^3$, we will consider, and still denote by $M^4(L^{(m)},d),$ the closed $4$-manifold obtained by adding a further $4$-handle.
 

\bigskip
Before describing the new procedure, the following preliminary  notations are needed:  
\begin{itemize}
\item[-] For each $i\in\{1,\ldots,m\},$ let us ``mark" two points $H_i$ and $H_i^{\prime}$ on $L_i$, such that they divide $L_i$ into two parts, one containing only overcrossings and the other containing only undercrossings of $L$.\footnote{Note that $1$-handles and $2$-handles may always be re-arranged, so to respect this requirement: see for example \cite[Prop. 4.2.7]{GS} or \cite[Chapter 1 - Principle 1]{[M]}.} 
\item[-] For each $j\in\{m+1,\ldots,l\},$ let us fix on $L_j$ a point $X_j$, between a curl and an undercrossing
 and let us consider the component $L_j$ in the diagram of $L$ as the union of consecutive {\it segments} obtained by cutting it not only at undercrossings, but also at overcrossings and at the point $X_j.$ 
\item[-] Then, for each $j \in \{m+1, \dots, l\}$ let us ``highlight" on $L_j$ - starting from $X_j$  and in the direction opposite to the undercrossing - a sequence $Y_j$ of consecutive segments, so that, for each  $i\in\{1,\ldots,m\},$  $H_i$ and $H_i^{\prime}$ belong to the boundary of the same region $\mathcal R_i$ of the ``diagram" obtained from $L$ by deleting the points $X_{m+1},\ldots, X_l$ and the segments of the sequences $Y_{m+1},\ldots, Y_l$ (with a little abuse of notation we will describe this new diagram as $L - \cup_{j=m+1}^l (X_j \cup Y_j)$).  
Note that $Y_j$ can be empty, while it never comprehends all segments of $L_j.$ 

Let us denote by $Y=Y_{m+1} \wedge \dots \wedge Y_l$ the sequence resulting from juxtaposition of the sequences of highlighted segments. 

\item[-] Finally, for each $i\in\{1,\ldots,m\}$,  let $\bar e_i$ (resp. $\bar e_i^{\prime}$) be the $1$-colored edge of $\Lambda(L,c)$ ``parallel" to the part of arc of $L_i$ containing the point $H_i$ (resp. $H_i^{\prime}$)  ``on the side" of the regions of $L$  merging into $\mathcal R_i$ (see Remark \ref{rem.drawn_over}), and let $v_i$ (resp.  $v_i^{\prime}$) be its endpoint belonging to the subgraph corresponding to an undercrossing of the dotted component $L_i$.
\end{itemize}

\bigskip

\medskip

\noindent 
{\bf  PROCEDURE C -  from $\mathbf{(L^{(m)},d)}$ to $\mathbf{\Gamma(L^{(m)},d)}$ \ (representing $\mathbf{M^4(L^{(m)},d)}$):}

\begin{itemize}
\item[(a)]  Let  $\Lambda(L,c)$  be the 4-colored graph constructed from $(L,c)$ according to Procedure $A$; in $\Lambda(L,c)$, let us choose a quadricolor $\mathcal Q_j$ for each (undotted)  component $L_j\ (j\in\{m+1,\ldots,l\})$ in the position corresponding to the point $X_j$. 
\item[(b)] Follow the sequence $Y=Y_{m+1} \wedge \dots \wedge Y_l$, starting, for each $j\in\{m+1,\ldots,l\}$ with $Y_j \ne \emptyset,$ from the segment corresponding to the pair of $1$-colored edges adjacent to vertices $P_4$ and $P_5$ identified by the quadricolor $Q_j;$ at each step of the sequence, if  $f,f^{\prime}$ is the pair of $1$-colored edges which are ``parallel" to the considered highlighted segment,  then:

if no 4-colored edge has already been added to the endpoints of $f$ and $f^\prime$, join, by $4$-colored edges, endpoints of $f$ to endpoints of $f^{\prime}$ belonging to different bipartition classes of $\Lambda(L,c);$ otherwise connect only the 
endpoints of $f$ and $f^{\prime}$ having no already incident $4$-colored edge. 

Moreover, if two consecutive highlighted segments correspond to an undercrossing, whose overcrossing does not correspond to previous segments in $Y$, add $4$-colored edges so to double the pairs of $0$-colored edges within the subgraph corresponding to that crossing.

\item[(c)] For each $i\in\{1,\ldots, m\}$,  add a 4-colored edge, so to connect $v_i$ and ${v_i^{\prime}}$.
\item[(d)] For each $j\in\{m+1,\ldots, l\}$, add a triad of $4$-colored edges between the vertices $P_{2r}$ and $P_{2r+1}$, $\forall r \in \{0,1,2\}$, of the quadricolor $\mathcal Q_j$ (as shown in Figure \ref{fig.quadricolor-singular}). 
\item[(e)] Add $4$-colored edges between the remaining vertices of $\Lambda(L,c)$, joining those which belong to the same $\{1,4\}$-residue.
\end{itemize}

\begin{rem}  {\em We point out that a quadricolor always arises in a component $L_j$ ($j \in \{m+1, \dots, l\}$) of $\Lambda(L,c)$ not only between a curl and an undercrossing but also between two curls with the same sign. Actually in this last case two quadricolors appear, one for each curl, and either of them can be indifferently chosen as $\mathcal Q_j$;
therefore we put the point $X_j$ between the curls and the sequence $Y_j$ can start from either ``side'' of it.
Moreover, note that the position of points $X_j$ may be suitably chosen, so to minimize the length of the sequence $Y$, 
provided that the above conditions for the existence of the quadricolor are satisfied.}
\end{rem}

\begin{figure}[!h]
\centerline{\scalebox{0.3}{\includegraphics{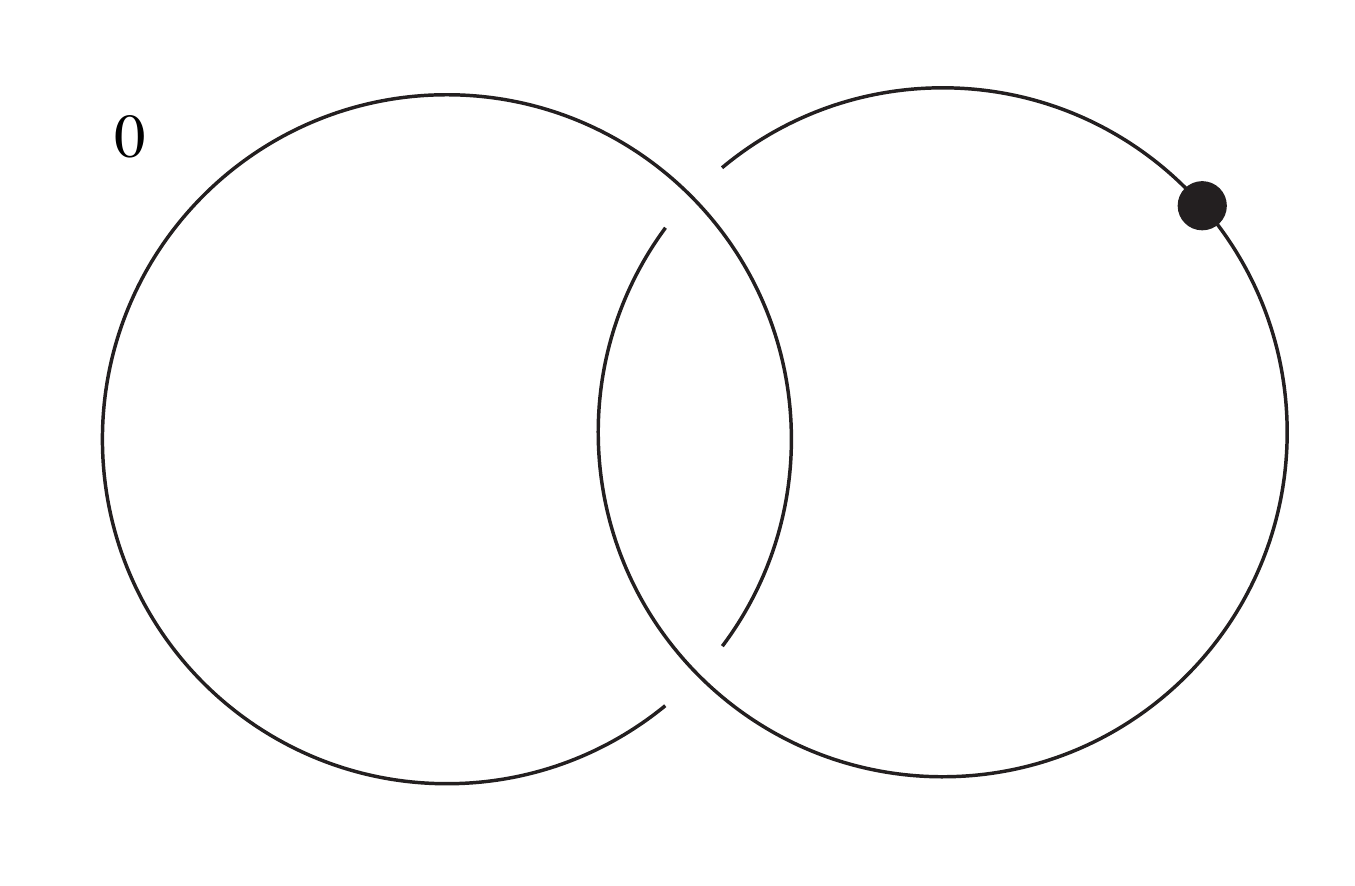}}}
\centerline{\scalebox{0.4}{\includegraphics{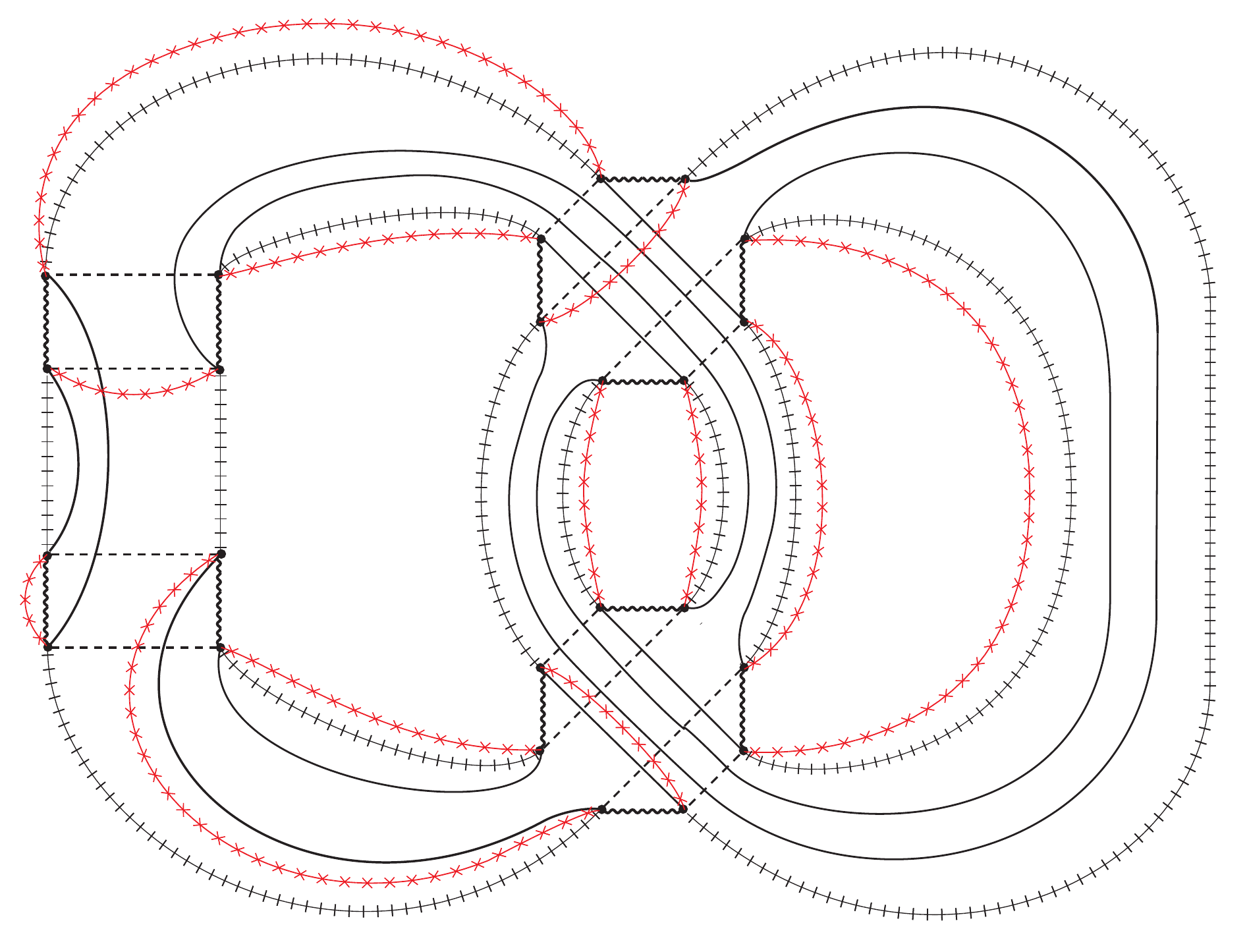}}}
\caption{{\footnotesize A Kirby diagram and the 5-colored graph representing the associated (closed) 4-manifold ($\mathbb S^4$)}}
\label{fig.dotted-Hopf}
\end{figure}

\begin{figure}[!h]
\centerline{\scalebox{0.4}{\includegraphics{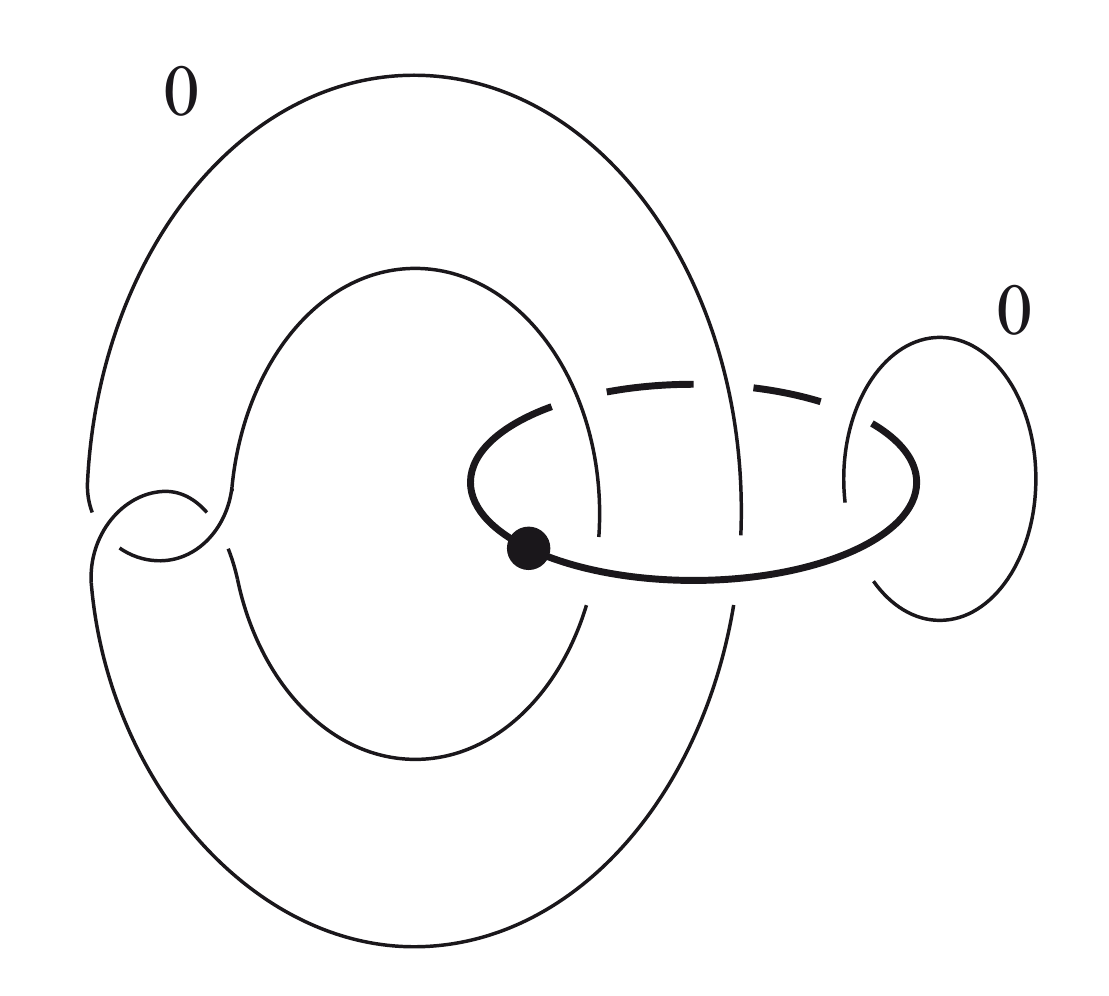}}}
\centerline{\scalebox{0.4}{\includegraphics{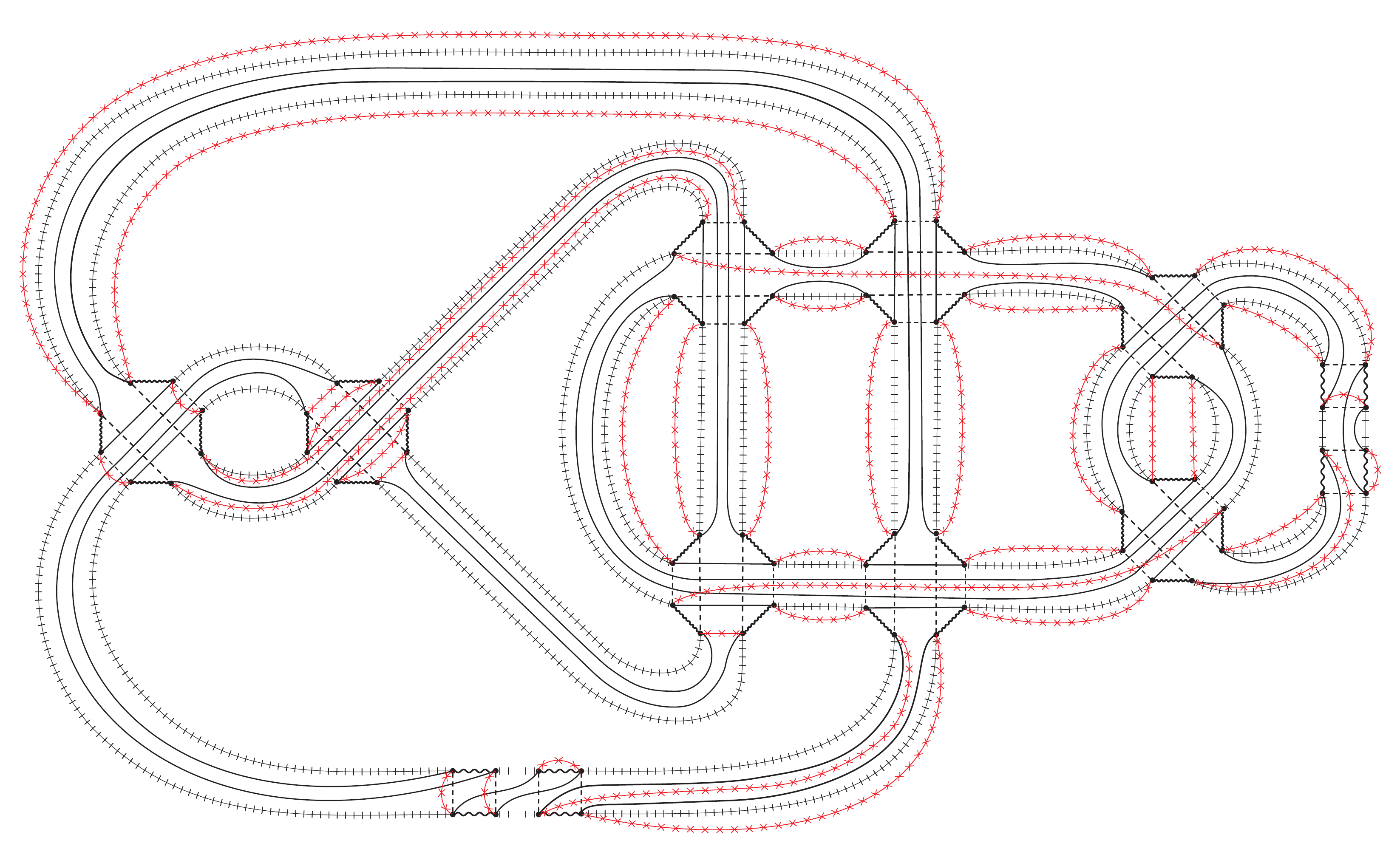}}}
\caption{{\footnotesize A Kirby diagram and the 5-colored graph representing the associated bounded 4-manifold ($\mathbb S^2 \times \mathbb D^2$)}}
\label{fig.fishtail}
\end{figure}

\begin{figure}[!h]
\centerline{\scalebox{0.45}{\includegraphics{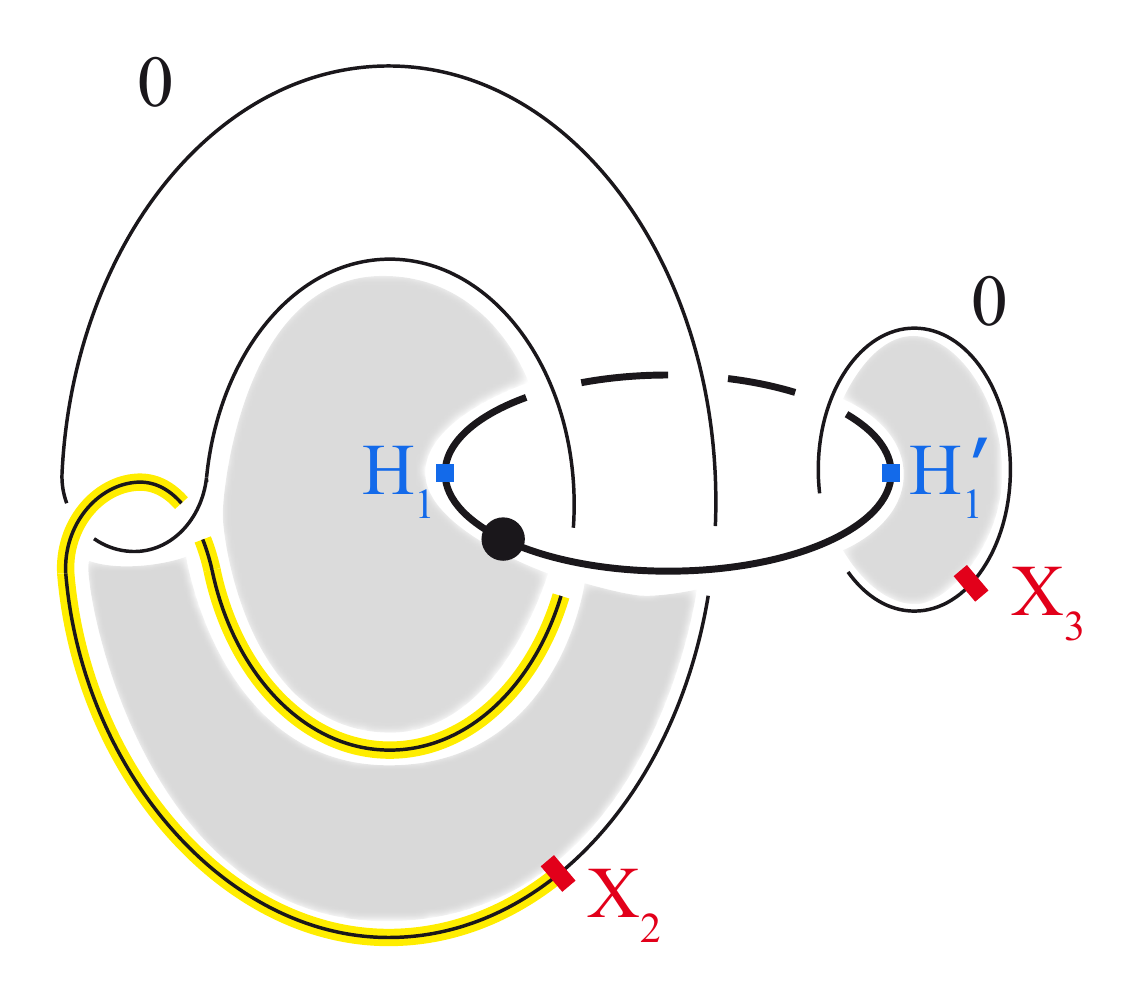}}}
\caption{{\footnotesize The  Kirby diagram of Fig. \ref{fig.fishtail}, with points and highlighted segments, according to Procedure C.  The yellow highlighted segments form the sequence $Y_2$, while $Y_3=\emptyset.$ The shaded regions, together with the infinite one, give rise to the region $\mathcal R_1$ of  $L - \cup_{j=2}^3 (X_j \cup Y_j)$ 
containing  both points $H_1$ and $H_1^{\prime}.$} }
\label{fig.fishtail_evidenziato}
\end{figure}

\begin{example} {\em Figures \ref{fig.dotted-Hopf} and \ref{fig.fishtail} show the result of the above construction applied to the depicted Kirby diagrams. In particular, note that step (b) of Procedure C is not required for the graph of Figure \ref{fig.dotted-Hopf},  since the highlighted sequence of segments is empty; on the contrary, the case of Figure \ref{fig.fishtail} requires to highlight a suitable set of consecutive segments in the Kirby diagram, as depicted in Figure  \ref{fig.fishtail_evidenziato}.
Via Kirby calculus, it is easy to check that the 5-colored graph in Figure \ref{fig.dotted-Hopf} represents the $4$-sphere, while the 5-colored graph in Figure \ref{fig.fishtail} represents $\mathbb S^2\times\mathbb D^2$; both facts can also be proved via suitable sequences of dipole moves. }
\end{example}

\begin{figure}[!h]
\centerline{\scalebox{0.47}{\includegraphics{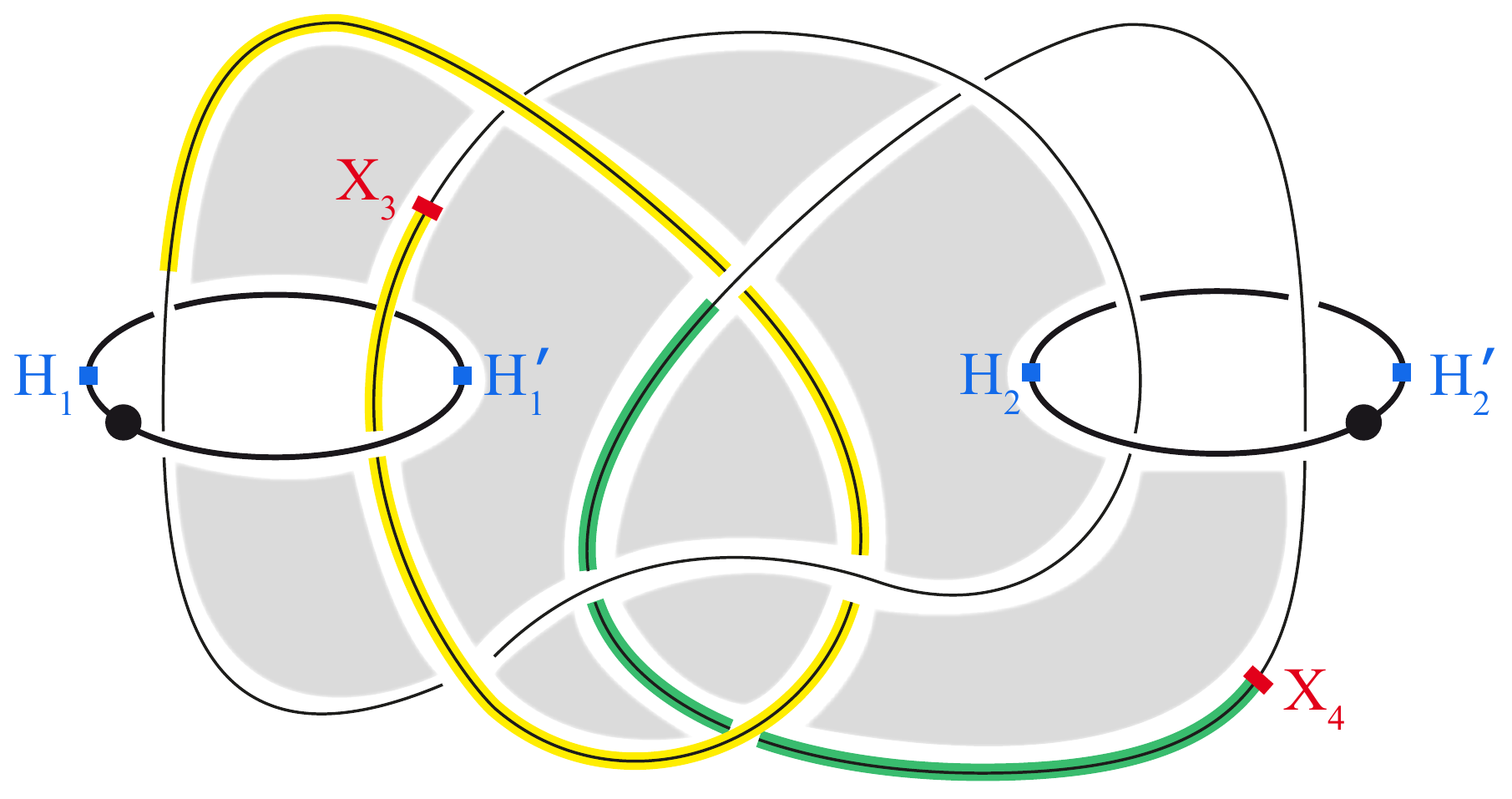}}}
\caption{}
\label{fig.K0&trefoil2}
\end{figure}

\begin{figure}[!h]
\centerline{\scalebox{0.47}{\includegraphics{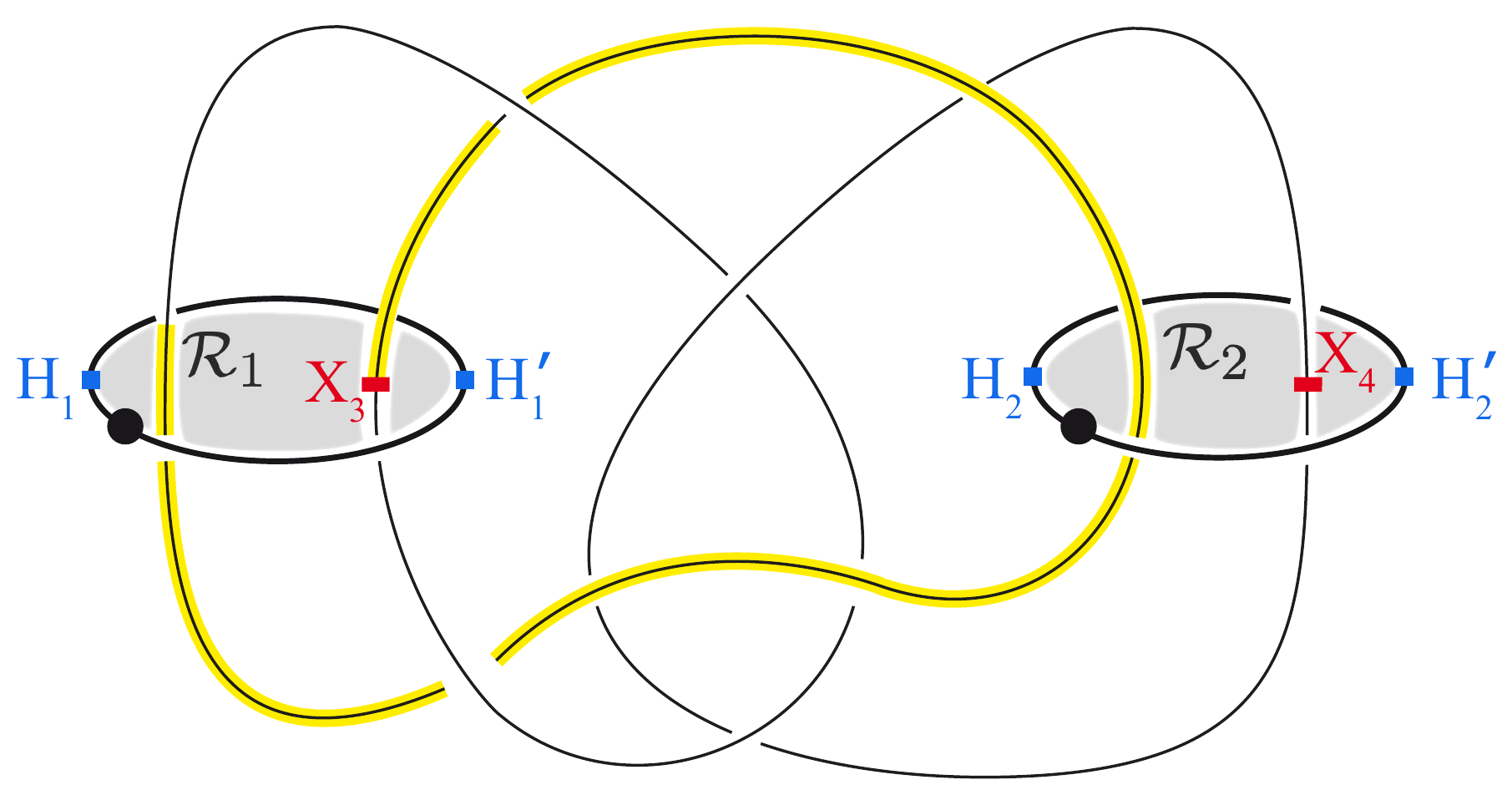}}}
\caption{}
\label{fig.K0&trefoil3}
\end{figure}

 \begin{example} {\em Given a framed link $(L^{(m)},d)$, the above construction may be implemented in different ways, depending on the choice of the points $X_i$ ($i=m+1,\dots, l$) on the framed components (step (a) of Procedure C) and on the choice of the sequence $Y=Y_{m+1} \wedge \dots \wedge Y_l$ of highlighted segments (step (b) of Procedure C). Figures \ref{fig.K0&trefoil2} and \ref{fig.K0&trefoil3} show two possibile ways to perform the above choices on the same Kirby diagram: in Figure \ref{fig.K0&trefoil2}  (resp. Figure \ref{fig.K0&trefoil3}) the yellow highlighted segments form the sequence $Y_3$,   while the green highlighted segments form the sequence $Y_4$ (resp. while $Y_4=\emptyset$). Note that, in the case of Figure \ref{fig.K0&trefoil2}, the regions $\mathcal R_1$ and $\mathcal R_2$ of $L - \cup_{j=3}^4 (X_j \cup Y_j)$ 
coincide: they are obtained by merging the shaded regions, together with the infinite one.  On the other hand, in the case of Figure \ref{fig.K0&trefoil3}, the regions $\mathcal R_1$ and $\mathcal R_2$ of $L - \cup_{j=3}^4 (X_j \cup Y_j)$ 
are distinct. }
\end{example}

The proof  that the graph obtained via Procedure C  really represents $M^4(L^{(m)},d)$ is given in Theorem \ref{thm.main(dotted)}. 
In order to help the reader, we can anticipate that it will be performed by means of the followings steps: 
\begin{itemize}
\item[(i)] starting from the $4$-colored graph $\Lambda(L,c)$ - already proved to represent $M^3(L,c)$ in \cite{Casali JKTR2000} - we obtain a $4$-colored graph $\Lambda_{smooth}$ representing $\#_m(\mathbb S^1 \times \mathbb S^2)$ by suitably exchanging a triad of $1$-colored edges for each framed component of $(L^{(m)},d)$ (Proposition \ref{Gamma_smooth_dotted}(i)); 
\item[(ii)] by capping-off with respect to color 1, we obtain a 5-colored graph representing  $[\#_m(\mathbb S^1 \times \mathbb S^2)]\times I$; 
\item[(iii)] this 5-colored graph is modified by a sequence of moves not affecting the represented 4-manifold (the so called $\rho_2$-pairs switching), so to have on one boundary component of $[\#_m(\mathbb S^1 \times \mathbb S^2)]\times I$ a particular structure  (called $\rho_3$-pair) for each dotted component of $(L^{(m)},d)$;
\item[(iv)] a suitable move ($\rho_3$-pair switching) is applied on each such structure, realizing - on the considered boundary component - the attachment of 1-handles corresponding to the $m$ dotted components of $(L^{(m)},d)$ (Proposition \ref{Gamma_smooth_dotted}(ii));
\item[(v)] by re-establishing  the triads of $1$-colored edges of step (i), the 5-colored graph $\Gamma(L^{(m)},d)$  is obtained. Since its only singular 4-residue is the $\hat 4$-residue $\Gamma(L,c)$, it represents a 4-manifold with connected boundary $M^3(L,c)$; moreover, $\Gamma(L^{(m)},d)$ represents  $M^4(L^{(m)},d)$  since - similarly as in Procedure B - each triad re-exchanging is proved to correspond to the addition of a 2-handle according to the framed component, on the remaining boundary component  (Proposition \ref{Gamma_smooth}(ii)). 
\end{itemize}

\medskip

In order to go into details, the notion of $\rho$-pair\footnote{$\rho$-pairs and their switching were introduced by Lins (\cite{Lins-book} and subsequently studied in \cite{Casali-Cristofori ElecJComb 2015}, \cite{Bandieri-Gagliardi}, \cite{CFMT}.} and some preliminary results are needed. 

\begin{defn}\label{rho-pair}  {\em A $\rho_h$-pair ($1\leq h\leq n$) of color $c\in\Delta_n$ in a bipartite $(n+1)$-colored graph $\Gamma$ is a pair of $c$-colored edges $(e,f)$ sharing the 
same $\{c,i\}$-colored cycle for each $i\in\{c_1,\ldots,c_h\}\subseteq\Delta_n.$ Colors $c_1,\ldots,c_h$ are said to be {\it involved}, while the other $n-h$ colors 
are said to be {\it not involved} in the $\rho_h$-pair.
\\
The {\it switching} of $(e,f)$ consists in canceling $e$ and $f$ and establishing new $c$-colored edges between their endpoints in such a way as to preserve the bipartition.}
\end{defn}

The topological consequences of the switching of $\rho_{n-1}$- and $\rho_n$-pairs have been completely determined in the case of closed $n$-manifolds: see \cite{Bandieri-Gagliardi}, where it is proved that a $\rho_{n-1}$-pair (resp. $\rho_n$-pair) switching does not affect the represented $n$-manifold  (resp. either induce the splitting into two connected summands, or the ``loss" of a $\mathbb S^1 \times \mathbb S^{n-1}$ summand in the represented $n$-manifold). In dimension three the study has been performed also in the case of manifolds with boundary: see \cite{CFMT}, where more cases are proved to occur.

As we will see in the proof of the following Proposition \ref{Gamma_smooth_dotted}, we are particularly interested in the effect of switching $\rho_{2}$- and 
$\rho_{3}$-pairs in $5$-colored graphs. A useful result is the following.

\begin{figure}[!h]
\centerline{\scalebox{0.55}{\includegraphics{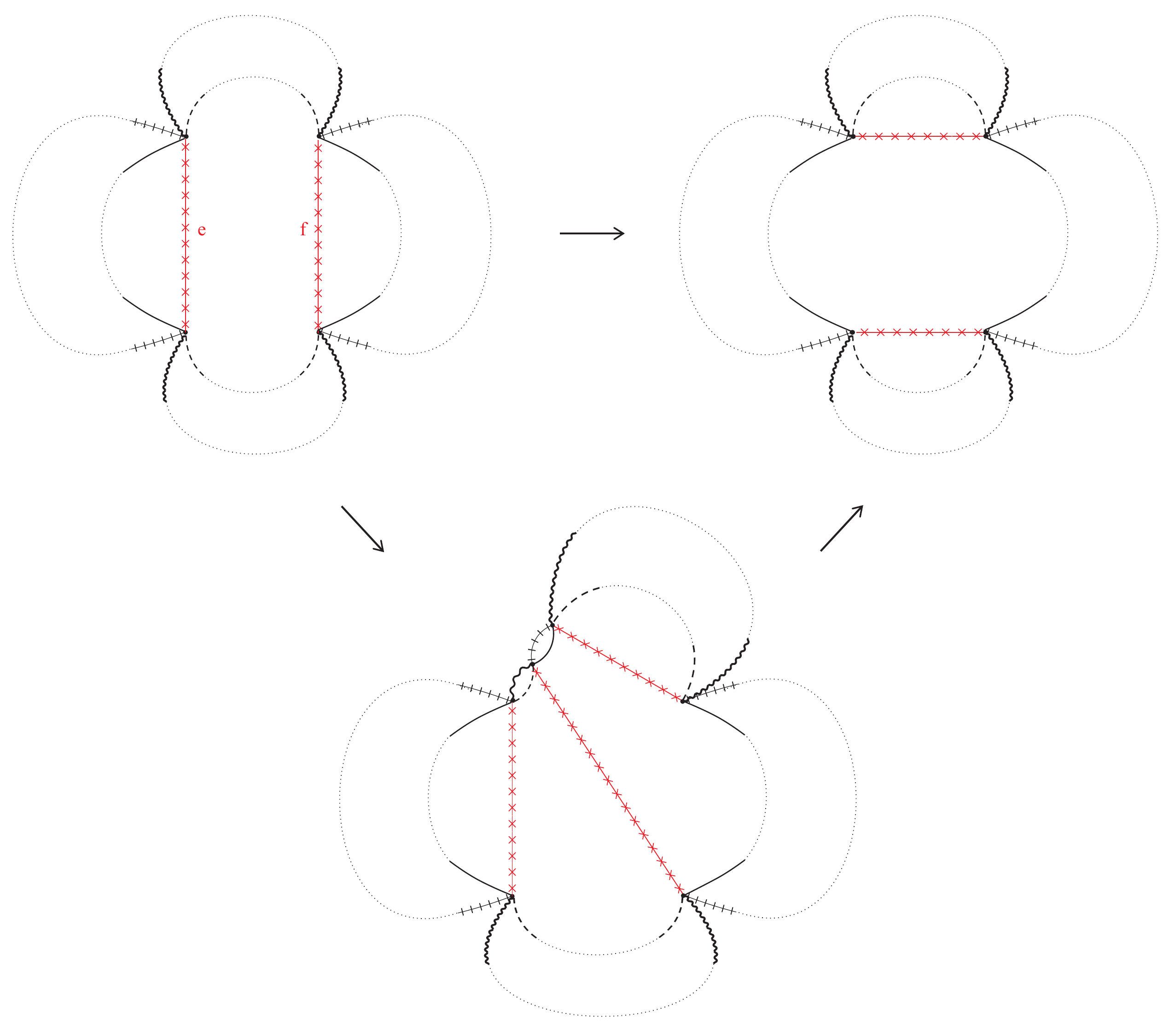}}}
\caption{{\footnotesize Factorization of a $\rho_2$-pair switching into two proper dipoles (not affecting the represented 4-manifold)} }
\label{fig: ro2} 
\end{figure}

\begin{lemma}\label{lemma-ro2-pairs} Let $(e,f)$ be a $\rho_2$-pair in a $5$-colored graph $\Gamma$ representing a compact $4$-manifold $M^4$ and let $\G^\prime$ be obtained from $\Gamma$ by switching the $\rho_2$-pair. Then $\Gamma^\prime$ represents $M^4$, too.

Moreover, for each cyclic permutation $\varepsilon$ of $\Delta_4$, where $\varepsilon_k$ is the color of $(e,f)$: 
\begin{itemize}
 \item [-] if both $\varepsilon_{k-1}$ and $\varepsilon_{k+1}$ are involved, then $\rho_\varepsilon(\Gamma^\prime)=\rho_\varepsilon(\Gamma)-1$
 \item [-] if neither $\varepsilon_{k-1}$ nor $\varepsilon_{k+1}$ is involved, then $\rho_\varepsilon(\Gamma^\prime)=\rho_\varepsilon(\Gamma)+1$
 \item [-] if exactly one between $\varepsilon_{k-1}$ and $\varepsilon_{k+1}$ is involved, then $\rho_\varepsilon(\Gamma^\prime)=\rho_\varepsilon(\Gamma).$
 \end{itemize}
\end{lemma}
\dimo In order to prove that $\Gamma^\prime$ represents $M^4$, too, it is sufficient to observe that the switching of $(e,f)$ can be factorized by a sequence of dipole moves as shown in Figure~\ref{fig: ro2}, i.e. by the addition of a $2$-dipole of the colors not involved in the $\rho_2$-pair,  followed by the cancellation of a $2$-dipole of the colors involved in the $\rho_2$-pair. Note that any $2$-dipole in a $5$-colored graph is proper (see Proposition \ref{proper-dipole}), and hence both moves do not change the represented manifold,  since - as already pointed out in Section \ref{prelim} - they correspond to re-triangulations of balls embedded in the cell-complexes associated to the involved colored graphs.

With regard to the regular genus of $\Gamma^\prime$ with respect to $\varepsilon$, note that the switching of $(e,f)$ increases by one (resp. decreases by one) the 
number of $\{\varepsilon_k,i\}$-colored cycles of $\Gamma$ if $i$ is an involved (resp. a not involved) color, while the number of $\{i,j\}$-colored cycles with 
$i,j\neq\varepsilon_k$ is not changed. An easy calculation yields the statement. 
\qed

\begin{prop}\label{Gamma_smooth_dotted} \  \par 
\begin{itemize}
\item[(i)]
The 4-colored graph $\Gamma^{(m)}_{smooth}$, obtained from $\Lambda(L,c)$ by 
exchanging the triad of $1$-colored edges (according to Figure \ref{fig.Gamma_smooth}) in a quadricolor for each framed component of $(L^{(m)},d)$, represents  $\#_m (\mathbb S^1 \times \mathbb S^2)$.  
\item[(ii)]
The 5-colored graph  $\bar\Gamma^{(m)}_{smooth}$,  obtained  by applying steps (b) and (c) of  Procedure C to $\Gamma^{(m)}_{smooth}$, and then by ``capping off" with respect to color $1$, represents  the genus $m$ 4-dimensional handlebody $\mathbb Y^4_m$.   
\end{itemize}
 \end{prop}

\dimo Part (i) directly follows from Proposition \ref{Gamma_smooth} (i), by noting that, if all framed components of $(L^{(m)},d)$ are deleted, only the $m$ dotted components remain, and the associated framed link, consisting in $m$ disjoint trivial 0-framed components, actually represents the 3-manifold $\#_m (\mathbb S^1 \times \mathbb S^2)$.

As regards part (ii), it is necessary to note that  $\Gamma^{(m)}_{smooth}$ gives rise, by ``capping off" with respect to color $1$, to a 5-colored graph representing $[\#_m (\mathbb S^1 \times \mathbb S^2)] \times I$, whose two boundary components - both homeomorphic to $\#_m (\mathbb S^1 \times \mathbb S^2)$ - are represented by the (color-isomorphic) 
subgraphs $\Theta$ and $\Theta^\prime$, obtained by deleting the $4-$colored and $1-$colored edges respectively. 
This 5-colored graph admits $\rho_2$-pairs of color $4$ in a suitable sequence induced by the sequence of $2$-dipoles whose cancellation from $\Gamma^{(m)}_{smooth}$ yields $\Lambda (\bigsqcup_m K_0, 0)$, the 4-colored graph associated to the trivial link with $m$ 0-framed components (see Remark \ref{rem.sequence_dipoles}, applied to all framed components of $(L^{(m)},d)$).
The  switching of these $\rho_2$-pairs is equivalent (up to ``capping off" with respect to color $1$) to the addition of  $4$-colored edges according to step (b) in $\Lambda(L,c).$

More precisely, the pairs of $4$-colored edges that have to be switched in the sequence of $\rho_2$-pairs are exactly the 4-colored edges adjacent to the pairs of vertices constituting  2-dipoles of the sequence of dipole eliminations  (starting, for each $j\in\{m+1,\ldots,l\}$ such that $Y_j\neq\emptyset$, with the dipole whose vertices are 2-adjacent to the vertices $P_4$ and $P_5$ identified by the quadricolor $\mathcal Q_j$) from  $\Gamma^{(m)}_{smooth}$ to $\Lambda (\bigsqcup_m K_0, 0)$; 
moreover, the colors involved in each $\rho_2$-pair are exactly those (never comprehending color $1$) of the corresponding $2$-dipole.  

Hence, the graph $\tilde \Gamma^{(m)}_{smooth}$, obtained  after all $\rho_2$-pairs switchings, still represents $[\#_m (\mathbb S^1 \times \mathbb S^2)] \times I$ and one of its boundary component is represented by  $\Theta$, too,  but the other is represented by the 4-colored graph  $\Theta^{\prime \prime}$ obtained from $\Theta^\prime$ by switching $\rho_2$-pairs induced by the above ones. 

We point out that, for each $i\in\{1,\ldots, m\}$, the pair of $4$-colored edges having an endpoint in $v_i$ and ${v^{\prime}_i}$ respectively, turn out to form a $\rho_3$-pair of color $4$ in $\tilde \Gamma^{(m)}_{smooth}.$
In fact, they double $\bar e_i$ and/or $\bar e_i^{\prime}$, or they arise from the possible switching of $4$-colored edges doubling $\bar e_i$ and/or $\bar e_i^{\prime}$ by the above sequence of $\rho_2$-pairs switchings; as a consequence, 
they belong both to the same $\{0,4\}$-residue and to the same $\{3,4\}$-residue (since  $\bar e_i$ and $\bar e_i^{\prime}$ share both the same $\{0,1\}$-residue and the same $\{1,3\}$-residue in $\Lambda(L,c)$), and the sequence of $\rho_2$-pair switchings makes them to belong also to the same $\{2,4\}$-residue (which corresponds  to the boundary of the region 
$\mathcal R_i$ of $L - \cup_{j=m+1}^l (X_j \cup Y_j)$). 

It is known that the switching of a $\rho_3$-pair in a 4-colored graph representing a closed 3-manifold has the effect of ``subtracting'' an $\mathbb S^1\times\mathbb S^2$ summand (see \cite{Bandieri-Gagliardi} for details); hence the switching of the above $m$ $\rho_3$-pairs transforms  $\Theta^{\prime \prime}$ into a 4-colored graph representing $\mathbb S^3$, while the $\hat 4$-residue $\Theta$ is unaltered and each $\hat c$-residue with $i\in\{0,2,3\}$ still represents the 3-sphere as in $\tilde \Gamma^{(m)}_{smooth}$ (since a $\rho_2$-pair switching has been performed in each affected $\hat c$-residue, for $i\in\{0,2,3\}$). 

\begin{figure}[!h]
\centerline{\scalebox{0.55}{\includegraphics{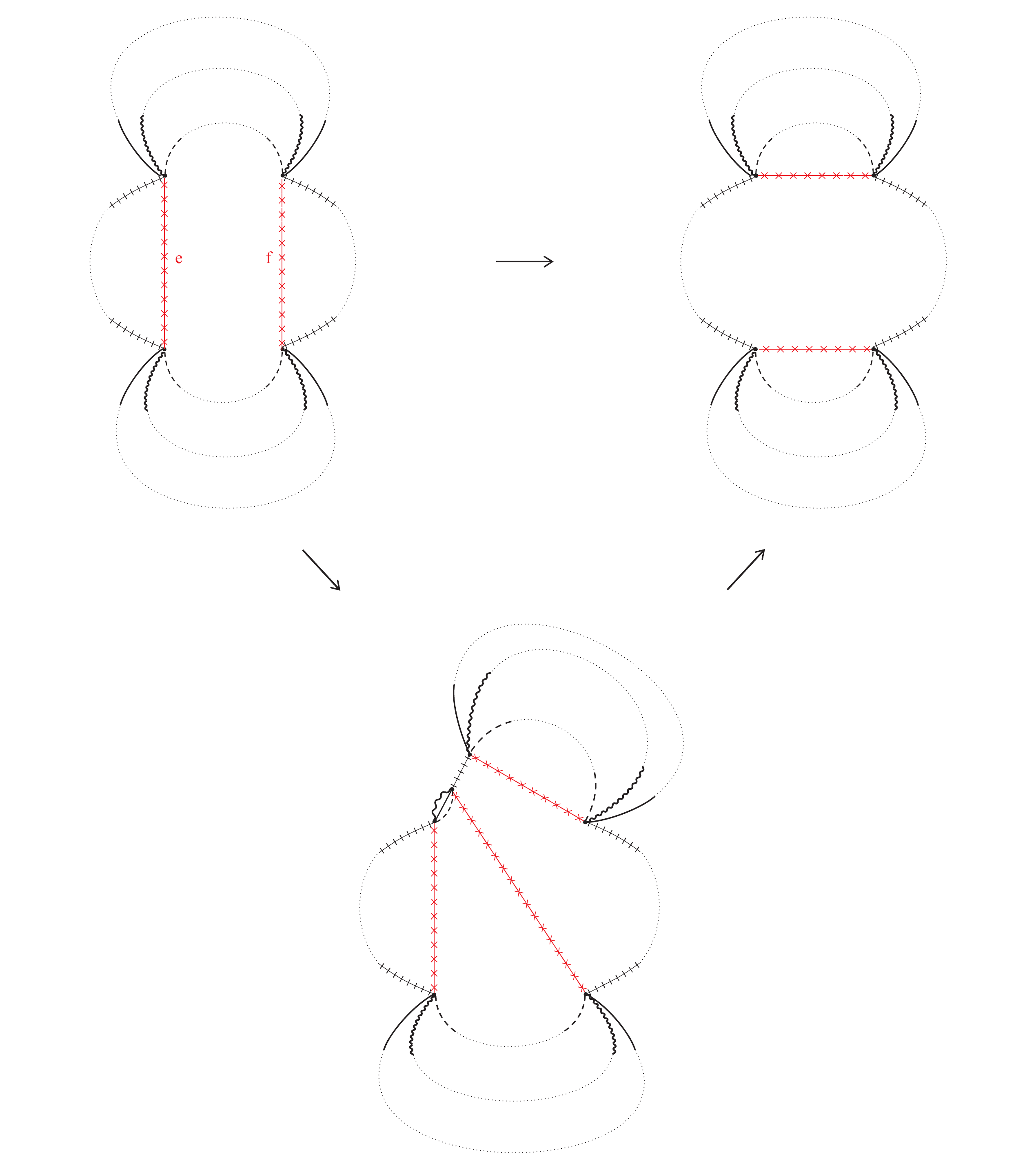}}}
\caption{{\footnotesize Factorization of a $\rho_3$-pair switching into two moves, the first (resp. second) one possibly affecting (resp. always not affecting) the represented 4-manifold}}
\label{fig: ro3} 
\end{figure}

Moreover, supposing $(e,f)$ to be one of the above $\rho_3$-pairs in $\tilde \Gamma^{(m)}_{smooth}$, its switching can be factorized as in Figure \ref{fig: ro3} by inserting a 1-colored edge and subsequently canceling a 3-dipole. 
The insertion of the $1$-colored edge in the colored triangulation associated to  $\tilde \Gamma^{(m)}_{smooth}$ consists in ``breaking" a tetrahedral face on the boundary of $[\#_m (\mathbb S^1 \times \mathbb S^2)] \times I$ corresponding to the $\hat 1$-residue $\Theta^{\prime \prime}$, and inserting a new pair of $4$-simplices sharing the same $3$-dimensional face opposite to the $1$-labelled vertex; hence, it  may be seen as the attachment of a polyhedron homeomorphic to $\mathbb D^3\times\mathbb D^1$ to the considered boundary,  so to transform it into a triangulation of $\#_{m-1} (\mathbb S^1 \times \mathbb S^2)$, without affecting the interior of $[\#_m (\mathbb S^1 \times \mathbb S^2)] \times I$, nor its boundary corresponding to the $\hat 4$-residue.\footnote{Actually, the switching of each of the $m$ $\rho_3$-pairs corresponds to the attaching of a $3$-handle to the boundary corresponding to the $\hat 4$-residue of $\tilde \Gamma^{(m)}_{smooth}$.}
Whenever all $m$ $\rho_3$-pairs are switched, the $\hat 1$-residue of the obtained $5$-colored graph comes to represent the $3$-sphere, i.e. the represented $4$-manifold has a connected boundary, corresponding to the (unaltered) $\hat 4$-residue   $\Theta=  \Gamma^{(m)}_{smooth}$. 

On the other hand, the  switching of these $\rho_3$-pairs is equivalent (up to ``capping off" with respect to color $1$) to the addition of $4$-colored edges  in $\Lambda(L,c)$ according to step (c); therefore, step (c) of  Procedure C can be thought as the identification $\phi$ between the boundary of a genus $m$ 4-dimensional handlebody $\mathbb Y^4_m$ and the 
boundary component represented by $\Theta^{\prime \prime}$ in the triangulation of $ [\#_m (\mathbb S^1 \times \mathbb S^2)] \times I$ obtained in step (b). 

This proves statement (ii), since  $\bar\Gamma^{(m)}_{smooth}$ - which admits $4$ as its unique singular color - turns out to represent    $\mathbb Y^4_m \cup_\phi ([\#_m (\mathbb S^1 \times \mathbb S^2)] \times I) \cong \mathbb Y^4_m.$
\qed
 
 \medskip

  \bigskip

We are now going to prove that the $5$-colored graph  $\G(L^{(m)},d)$ obtained via Procedure C (i.e. by applying to $\Lambda(L,c)$ steps (b)-(e)) represents the compact 4-manifold associated to the Kirby diagram;  we will also give an estimation of its regular genus and compute its order. With this aim, if $(L^{(m)},d)$ is a Kirby diagram with $l$ components where the first $m > 0$ ones are dotted and $s$ crossings,  and $(L,c)$ is its associated framed link, let us set, for each $i\in\{m+1,\ldots,l\}$, 
$$\bar t_i = \begin{cases}
                                                       \vert w_i-c_i\vert\quad  if\ w_i \ne c_i  \\  
                                                      2  \quad otherwise
                                                      \end{cases}$$ where $w_i$ denotes the writhe of the $i$-th (framed) component of $L;$
moreover, let us denote by $u$ the number of undercrossings which are passed when following the sequence $Y$,  with the condition that the associated overcrossing does not correspond to 
previous segments in the sequence itself.

\medskip

\begin{thm}\label{thm.main(dotted)}
For each Kirby diagram $(L^{(m)},d)$, 
the bipartite $5$-colored graph $\G(L^{(m)},d)$ represents the compact $4$-manifold $M^4(L^{(m)},d).$ 

Moreover, it has regular genus  less or equal to 
\ $s + (l-m) + u +1$ and, if $L$ is different from the trivial knot, its order is \ $8s + 4\sum_{i=m+1}^l \bar t_i.$ 
\end{thm}

\dimo
In order to prove the first statement, we point out that in the proof of Proposition \ref{Gamma_smooth_dotted} (ii) we have considered a suitable triangulation of  $[\#_m (\mathbb S^1 \times \mathbb S^2)] \times I$, and then we have ``closed" one of its boundary components by identifying it with the boundary of the genus $m$ 4-dimensional handlebody  (via the addition of 4-colored edges according to steps (b) and (c) of the Procedure C). Hence, the polyhedron represented by $\bar\Gamma^{(m)}_{smooth}$ may be seen as the union of $0$- and $1$-handles of  $M^4(L^{(m)},d),$ plus a ``collar" on its boundary . 
Moreover, the ``free" boundary, homeomorphic to $\#_m (\mathbb S^1 \times \mathbb S^2)$, is represented by the 4-colored graph $\Theta= (\bar\Gamma^{(m)}_{smooth})_{\hat 4}= \Gamma^{(m)}_{smooth}$.  
Then, in order to obtain a $5$-colored graph representing $M^4(L^{(m)},d),$ it is sufficient to operate on this ``free" boundary, so perform the addition of a $2$-handle according to each framed component of $(L^{(m)},d)$. 
Now, the proof of Proposition \ref {Gamma_smooth}(ii) shows that the goal is achieved by exchanging the triad of $1$-colored edges, according to Figure \ref{fig.2-handle}, in the quadricolor $Q_j$ of the $j$-th component of $(L^{(m)},d)$, for each $j\in\{m+1,\ldots, l\}$. 
Since all these exchanging of $1$-colored edges have the effect to transform $\Gamma^{(m)}_{smooth}$ into $\Lambda(L,c)$, and step (d) applied to $\Lambda(L,c)$ is equivalent to the exchanging of $1$-colored edges according to Figure \ref{fig.2-handle} applied to  $\bar\Gamma^{(m)}_{smooth}$,  the final 5-colored graph representing the compact $4$-manifold $M^4(L^{(m)},d)$ turns out to be obtained by applying directly to $\Lambda(L,c)$ steps (b)-(d),  
and then  by ``capping off" with respect to color $1$ (step (e)).  

\medskip
In order to give an estimation of  the regular genus of $\G(L^{(m)},d)$, we first recall that $\rho_{\bar\varepsilon}(\Lambda(L,c))=s+1$ with $\bar\varepsilon=(1,0,2,3)$ (see also the proof of Theorem \ref{M4(L,c)}(i)), and hence that $s+1$ is also the regular genus, with respect to the permutation $\varepsilon = (1,0,2,3,4)$, of the $5$-colored graph obtained by doubling the 1-colored edges of $\Lambda(L,c)$ by color $4$.
Then, we have to analyze how the the regular genus is affected by the switchings of $\rho_2$- and $\rho_3$-pairs and the exchanging of triad of edges in the quadricolors described in the proofs of Proposition \ref{Gamma_smooth_dotted} and Theorem \ref{M4(L,c)}.

Now, let us point out that color $1$ is never involved in the considered $\rho_2$-pairs, while color $3$ is involved only in one of the two $\rho_2$-pairs corresponding to an undercrossing whose associated overcrossing does not correspond to previous segments in the sequence $Y$. Therefore, by Lemma \ref{lemma-ro2-pairs}, the regular genus with respect to $\varepsilon$ increases by $u$, when performing the sequence of $\rho_2$-pairs corresponding to the sequence $Y.$

With regard to the $\rho_3$-pairs, since they do not involve color $1$, which is consecutive in $\varepsilon$ to color $4$, the same argument used in the proof of Lemma \ref{lemma-ro2-pairs}, shows that the regular genus does not change after their switchings.

Finally, the exchanging of the triad of $4$-colored edges in a quadricolor, producing the attaching of a $2$-handle, decreases by two the number of $\{1,4\}$-colored cycles, while the numbers of all other bicolored cycles remain unaltered (see Figure \ref{fig.2-handle}). Hence, the regular genus increases by one for each quadricolor. Since the quadricolors are $l-m$, the statement is proved. 

\medskip

The proof of the theorem is completed  by noting that $\G(L^{(m)},d)$  has exactly the same order as  $\G(L,c)$ (and as $\Lambda(L,c)$, too). 
Hence, its calculation directly follows from  Theorem \ref{M4(L,c)} (i). 
\qed

\bigskip

We are now able to prove both upper bounds for the invariants of the 4-manifold associated to a Kirby diagram, already stated in Theorem \ref{regular-genus&gem-complexity(dotted)} in the Introduction.

\bigskip

\noindent\textit{Proof of Theorem~{\upshape\ref{regular-genus&gem-complexity(dotted)}}}
The upper bound for the regular genus of $M^4(L^{(m)},d)$ directly follows from the computation of $\rho_{\bar\varepsilon}(\G(L^{(m)},d))$ obtained in Theorem {\ref{thm.main(dotted)}}, together with the trivial inequality $u \le \bar s$. 

As regards the upper bound for the gem-complexity, it is sufficient to make use of the computation of the order of $\G(L^{(m)},d)$ obtained in Theorem {\ref{thm.main(dotted)}}, by pointing out that $\G(L^{(m)},d)$  contains a pair of $3$-dipoles of colors $\{0,1,4\}$ for each pair of adjacent undercrossings of dotted components; hence, a new $5$-colored graph $\G^{\prime}(L^{(m)},d))$ representing  $M^4(L^{(m)},d)$ may be obtained, with $$ \#V(\G^{\prime}(L^{(m)},d)) =  \#V(\G(L^{(m)},d)) - 4 [ (s- \bar s) - m]  =  4s + 4 \bar s + 4m + 4\sum_{i=m+1}^l \bar t_i.$$\qed

\bigskip
\bigskip

\noindent\textbf{Acknowledgements.}  This work was supported by GNSAGA of INDAM and by the University of Modena and Reggio Emilia, project: {\it ``Discrete Methods in Combinatorial Geometry and Geometric Topology".     
}

\end{document}